\documentclass[a4paper,11pt]{amsart}
\usepackage{amsmath,amsthm,amssymb}
\usepackage[mathscr]{eucal}
 \usepackage{cite}
\usepackage{upgreek}
\usepackage[bookmarks=false]{hyperref}
\usepackage{enumerate}
\usepackage{amsmath}
\usepackage{mathrsfs}
\usepackage{blindtext}
\usepackage{scrextend}
\usepackage{enumitem}
\addtokomafont{labelinglabel}{\sffamily}
\usepackage{color}
\usepackage[at]{easylist}

\numberwithin{equation}{section}

\theoremstyle{plain}
\newtheorem{main}{Theorem}
\newtheorem{mcor}[main]{Corollary}

\newtheorem{theorem}{Theorem}[section]
\newtheorem{claim}[theorem]{Claim}
\newtheorem{lemma}[theorem]{Lemma}
\newtheorem{proposition}[theorem]{Proposition}
\newtheorem{corollary}[theorem]{Corollary}
\theoremstyle{definition}
\newtheorem{definition}[theorem]{Definition}
\newtheorem*{definition*}{Definition}

\newtheorem{notation}[theorem]{Notation}
\newtheorem{construction}[theorem]{Construction}
\newtheorem{remark}[theorem]{Remark}
\newtheorem{fact}[theorem]{Fact}

\begin{document}

\title[Almost commuting matrices and stability for product groups]
{Almost commuting matrices and \\ stability for product groups}

\author[A. Ioana]{Adrian Ioana}
\address{Department of Mathematics, University of California San Diego, 9500 Gilman Drive, La Jolla, CA 92093, USA}
\email{aioana@ucsd.edu}

{\thanks{The author was supported in part by 
NSF FRG Grant \#1854074.}}
\begin{abstract} 
We prove that any product of two non-abelian free groups, $\Gamma=\mathbb F_m\times\mathbb F_k$, for $m,k\geq 2$,
 is not Hilbert-Schmidt stable. This means that there exist asymptotic representations $\pi_n:\Gamma\rightarrow \text{U}({d_n})$ with respect to the normalized Hilbert-Schmidt norm 
 which are not close to actual representations.
As a consequence, we prove the existence of contraction matrices $A,B$  such that $A$ almost commutes with $B$ and $B^*$, with respect to the normalized Hilbert-Schmidt norm, but $A,B$ are not close to any matrices $A',B'$ such that $A'$ commutes with $B'$ and $B'^*$.  This settles in the negative a natural version of a question concerning almost commuting matrices 
posed by Rosenthal in 1969. 
 \end{abstract}

\maketitle

\section{Introduction and statement of main results}
A famous question, which can be traced back to the foundations of quantum mechanics \cite{vN29}, is whether two matrices $A,B$,
which almost commute with respect to a given norm, 
must be close to two commuting matrices $A',B'$.  It was first explicitly posed by Rosenthal \cite{Ro69} for the normalized Hilbert-Schmidt norm and by Halmos \cite{Ha76} for the operator norm.  
Almost commuting matrices have since been studied extensively and found applications to several areas of mathematics, including operator algebras and group theory, quantum physics and computer science (see, e.g, the introductions of \cite{LS13,ES19}).
The most interesting case of this question is when the matrices 
 are contractions, and ``almost" and ``close" are taken independent of their sizes. 
The answer depends both on the types of matrices considered and the norms chosen. Historically, research has focused on the operator norm. In this situation, the answer is positive for self-adjoint matrices by a remarkable result of Lin \cite{Li97} (see also \cite{FR96,Ha08,KS14}), but negative for unitary \cite{Vo83} and general matrices \cite{Ch88} (see \cite{Da85,EL89} for related results). 
More recently, several works  \cite{HL08,Gl10,FK10,FS11,Sa14,HS16,HS17} studied the question for the normalized Hilbert-Schmidt norm and obtained affirmative answers for pairs of self-adjoint, unitary and normal matrices. In fact, the answer is positive if at least one of the matrices is normal, see Remark \ref{stable}(1). However, these results leave wide open the general situation when neither matrix is normal.

We make progress on this problem by proving that, in contrast to the case of normal matrices, a version of Rosenthal's question \cite{Ro69} has a negative answer for non-normal matrices.   The version that we consider is natural from the perspective of (self-adjoint) operator algebras.  Indeed, it requires that $A$ almost commutes not only with $B$ but also with its adjoint, $B^*$.

\begin{main}\label{acm}
There exist sequences of matrices $A_n,B_n\in\mathbb M_{d_n}(\mathbb C)$, for some $d_n\in\mathbb N$, such that 
\begin{enumerate}[label=(\alph*)]
\item $\|A_n\|,\|B_n\|\leq 1$, for every $n\in\mathbb N$, 
\item $\lim\limits_{n\rightarrow\infty}\|A_nB_n-B_nA_n\|_2=\lim\limits_{n\rightarrow\infty}\|A_nB_n^*-B_n^*A_n\|_2=0$ and
\item $\inf_{n\in\mathbb N}\big(\|A_n-A_n'\|_2+\|B_n-B_n'\|_2)>0$, for any sequences of matrices $A_n',B_n'\in\mathbb M_{d_n}(\mathbb C)$ such that $A_n'B_n'=B_n'A_n'$ and $A_n'B_n'^*=B_n'^*A_n'$, for every $n\in\mathbb N$.
\end{enumerate}
\end{main}

For $A={(a_{i,j})}_{i,j=1}^n\in\mathbb M_n(\mathbb C)$, we denote by $\|A\|$, $\|A\|_2=\big(\frac{1}{n}\sum_{i,j=1}^n|a_{i,j}|^2\big)^{\frac{1}{2}}$ and $\tau(A)=\frac{1}{n}\sum_{i=1}^na_{i,i}$ the operator norm,  normalized Hilbert-Schmidt norm and normalized trace of $A$.

\begin{remark}\label{stable}  We continue with two remarks on the statement of Theorem \ref{acm}. 
\begin{enumerate}
\item The conclusion of Theorem \ref{acm} fails if one of the matrices is normal. 
Moreover,  the following holds: let $A_n,B_n\in\mathbb M_{d_n}(\mathbb C)$  be contractions such that $\|A_nB_n-B_nA_n\|_2\rightarrow 0$ and $B_n$ is normal, for every $n\in\mathbb N$. Then there are $A_n',B_n'\in\mathbb M_{d_n}(\mathbb C)$ such that $A_n'B_n'=B_n'A_n'$ and $A_n'B_n'^*=B_n'^*A_n'$, for every $n\in\mathbb N$, and   $\|A_n-A_n'\|_2+\|B_n-B_n'\|_2\rightarrow 0$
(see Lemma \ref{ultra}).

\item Theorem \ref{acm} complements a result of von Neumann \cite[Theorem 9.7]{vN42} 
which implies the existence of contractions $A_n\in\mathbb M_{k_n}(\mathbb C)$, for some $k_n\rightarrow\infty$, such that any contractions $B_n\in\mathbb M_{k_n}(\mathbb C)$ which verify condition (b),  must satisfy that $\|B_n-\tau(B_n)1\|_2\rightarrow 0$. In particular, $A_n,B_n$ are close to the commuting matrices $A_n,\tau(B_n)1$. Thus, the pair $A_n,B_n$ does not satisfy the conclusion of Theorem \ref{acm}, for any choice of contractions $B_n\in\mathbb M_{k_n}(\mathbb C)$. Moreover, the arguments used in \cite{vN42} 
are probabilistic, 
which suggests
that randomly chosen contractions $A_n\in\mathbb M_{k_n}(\mathbb C)$ should have the above property. 
This rules out examples as in Theorem \ref{acm}, where one of the matrices is chosen randomly.
Nevertheless, we use matrices satisfying (a strengthening of) the property of \cite{vN42} as building blocks in our construction of $A_n,B_n$  as in Theorem \ref{acm} (see the comments at the end of the introduction).

\end{enumerate}
\end{remark}

Theorem \ref{acm} is a consequence of a non-stability result for the product group $\mathbb F_2\times\mathbb F_2$, see Theorem \ref{main}. To motivate the latter result,
 we note that whether almost commuting matrices are near commuting ones is a prototypical stability problem. In general, following \cite{Hy41,Ul60}, stability refers to a situation when elements which ``almost" satisfy an equation must be ``close" to elements satisfying the equation exactly.
  In recent years, there has been a considerable amount of interest in the study of group stability  (see \cite{Th17,Io19}). 
For a countable group $\Gamma$, one can define stability with respect to any class $\mathcal C$ of metric groups endowed with bi-invariant metrics. This requires that any asymptotic homomorphism from $\Gamma$ to a group in $\mathcal C$ is close to an actual homomorphism  \cite{AP14,AP17,CGLT17,Th17}. Specializing to the class $\mathcal C$ of unitary groups  endowed with the normalized Hilbert-Schmidt norms leads to the following notion of stability introduced in \cite{HS17,BL18}: 

\begin{definition}
A sequence of maps $\varphi_n:\Gamma\rightarrow\text{U}(d_n)$, for some $d_n\in\mathbb N$, is called an {\it asymptotic homomorphism} if it satisfies $\lim\limits_{n\rightarrow\infty}\|\varphi_n(gh)-\varphi_n(g)\varphi_n(h)\|_2=0$, for every $g,h\in\Gamma$. The group $\Gamma$ is called {\it Hilbert-Schmidt stable} (or \text{HS}-{\it stable}) if for any asymptotic homomorphism $\varphi_n:\Gamma\rightarrow\text{U}(d_n)$, we can find homomorphisms $\rho_n:\Gamma\rightarrow\text{U}(d_n)$ such that $\lim\limits_{n\rightarrow\infty}\|\varphi_n(g)-\rho_n(g)\|_2=0$, for every $g\in\Gamma$. 
\end{definition}

 The class of HS-stable groups includes the free groups $\mathbb F_m$, virtually abelian groups and one-relator groups with non-trivial center  \cite{HS17}, certain graph product groups \cite{At18}, and is closed under free products. 
Moreover, the product of two HS-stable groups is HS-stable, provided that one of the groups is abelian \cite[Theorem 1]{HS17} or, more generally, amenable \cite[Corollary D]{IS19}. 

However, it remained a basic open problem whether HS-stability is closed under general direct products and, specifically, if $\mathbb F_2\times\mathbb F_2$ is HS-stable (see \cite[Remark 1.4]{Io19}).
We settle this problem by proving that the product of two non-abelian free groups is not HS-stable.  Moreover, we show:

\begin{main}\label{main}
The group $\mathbb F_k\times\mathbb F_m$ is not flexibly \emph{HS}-stable, for any integers $k,m\geq 2$.
\end{main}

Before discussing  the notion of flexible HS-stability and results related to Theorem \ref{main}, let us outline how Theorem \ref{main} implies Theorem \ref{acm}. Let $\varphi_n:\mathbb F_2\times\mathbb F_2\rightarrow\text{U}(d_n)$ be an asymptotic homomorphism which witnesses that $\mathbb F_2\times\mathbb F_2$ is not HS-stable and $a_1,a_2\in\mathbb F_2$ be free generators. For $1\leq j\leq 2$ and $n\in\mathbb N$, let $h_{n,j},k_{n,j}\in\mathbb M_{d_n}(\mathbb C)$ be self-adjoint matrices  with spectrum contained in $[-\frac{1}{2},\frac{1}{2}]$ such that $\varphi_n(a_j,e)=\exp(2\pi ih_{n,j})$ and $\varphi_n(e,a_j)=\exp(2\pi i k_{n,j})$. We then prove that the matrices $A_n=h_{n,1}+ik_{n,1}$ and $B_n=h_{n,2}+ik_{n,2}$ satisfy the conclusion of Theorem \ref{acm}.

It was shown in \cite{BL18} that all infinite residually finite property (T) groups $\Gamma$ (e.g., $\text{SL}_n(\mathbb Z)$, for $n\geq 3$), are not HS-stable.
The proof  builds on the observation that any sequence of homomorphisms $\rho_n:\Gamma\rightarrow\text{U}(d_n)$ with $d_n\rightarrow\infty$ can be perturbed slightly to obtain an asymptotic homomorphism $\varphi_n:\Gamma\rightarrow\text{U}(d_n-1)$. To account for this method of constructing asymptotic homomorphisms, the following weakening of the notion of HS-stability was suggested in \cite{BL18}:

\begin{definition}
A countable group $\Gamma$ is called {\it flexibly} HS-{\it stable} if for any asymptotic homomorphism $\varphi_n:\Gamma\rightarrow\text{U}(d_n)$, we can find homomorphisms $\rho_n:\Gamma\rightarrow\text{U}(D_n)$, for some $D_n\geq d_n$ such that $\lim\limits_{n\rightarrow\infty}\frac{D_n}{d_n}=1$ and  $\lim\limits_{n\rightarrow\infty}\|\varphi_n(g)-p_n\rho_n(g)p_n\|_2=0$, for every $g\in\Gamma$, where $p_n:\mathbb C^{D_n}\rightarrow\mathbb C^{d_n}$ denotes the orthogonal projection for every $n\in\mathbb N$.
\end{definition}

If a Connes-embeddable countable group $\Gamma$ is flexibly HS-stable, then it must be residually finite. On the other hand, deciding if a residually finite group is flexibly HS-stable or not is a challenging problem. For instance, while the arithmetic groups $\text{SL}_n(\mathbb Z), n\geq 3,$ are not HS-stable
 by \cite{BL18}, it is open whether they are flexibly HS-stable. The first examples of residually finite groups which are not flexibly HS-stable were found only recently in \cite{ISW20}, where certain groups with  the relative property (T), including $\mathbb Z^2\rtimes\text{SL}_2(\mathbb Z)$, were shown to have this property.
 
 Theorem \ref{main} provides the only other known examples of non-flexibly HS-stable residually finite groups, and the first  that do not have infinite subgroups with the relative property (T).
Moreover, these are the first examples of residually finite non-HS-stable groups that neither satisfy property (T;FD) (see \cite[Section 4.2]{BL18}) nor have infinite subgroups with the relative property (T).

\begin{remark} We now compare Theorem \ref{main} with two related results concerning other notions of stability.
A countable group  $\Gamma$ is called \text{P}-{\it stable} if it is stable with respect to the class of finite permutation groups endowed with the normalized Hamming distance (see \cite{Io19} for a survey on P-stability). 
As shown in \cite[Corollary B]{Io19}, P-stability is not closed under direct products. Given the similarity between the notions of HS-stability and P-stability \cite{AP14}, it should not be surprising that HS-stability is not closed under direct products. 
We note however that the methods of \cite{Io19} cannot be adapted to prove Theorem \ref{main}. 
The approach of \cite{Io19}, which exploits the discrete aspects of P-stability, allows to prove that the group $\mathbb F_2\times\mathbb Z$ is not P-stable, despite being HS-stable by \cite[Theorem 1]{HS17}. 
As we explain at the end of the introduction, to prove Theorem \ref{main} we introduce an entirely new approach based on ideas from the theory of von Neumann algebras.

A countable group $\Gamma$ is called {\it W$^*$-tracially stable} if it is stable with respect to the class of unitary groups of tracial von Neumann algebras endowed with their $2$-norms \cite{HS17}.   Theorem \ref{main}  strengthens  \cite[Theorem E]{IS19} which showed $\mathbb F_k\times\mathbb F_m$ is not W$^*$-tracially stable, for any integers $k,m\geq 2$. Indeed, being W$^*$-tracial stable is  stronger than being HS-stable, which corresponds to restricting to unitary groups of finite dimensional von Neumann algebras.

\end{remark} 

Next, we mention two reformulations of Theorem \ref{main} in terms of operator algebras. Let $(M_n,\tau_n)$, $n\in\mathbb N$, be a sequence of tracial von Neumann algebras and $\omega$ be a free ultrafilter on $\mathbb N$. The {\it tracial ultraproduct von Neumann algebra} $\prod_\omega M_n$ is defined as the quotient $\ell^{\infty}(\mathbb N,M_n)/\mathcal I_\omega(\mathbb N,M_n)$ of the C$^*$-algebra $\ell^{\infty}(\mathbb N,M_n)$ of sequences $(x_n)\in \prod_\mathbb N M_n$ with $\sup\|x_n\|<\infty$ by its ideal $\mathcal I_\omega(\mathbb N,M_n)$ of sequences $(x_n)$ such that $\lim\limits_{n\rightarrow\omega}\|x_n\|_2=0$.

First, by \cite[Proposition C]{IS19}, if $P, Q$ are commuting separable subalgebras of a tracial ultraproduct $\prod_\omega M_n$, and $P$ is amenable, then there are commuting von Neumann subalgebras $P_n,Q_n$ of $M_n$, for all $n\in\mathbb N$, such that $P\subset\prod_\omega P_n$ and $Q\subset\prod_\omega Q_n$.
In contrast, Theorem \ref{main} implies that, without the amenability assumption, this lifting property fails in certain matricial ultraproducts: 

\begin{mcor}\label{C}
There exist a sequence $(d_n)\subset\mathbb N$ and commuting separable von Neumann subalgebras $P,Q$ of $\prod_\omega\mathbb M_{d_n}(\mathbb C)$ such that the following holds: there are no commuting von Neumann subalgebras $P_n,Q_n$ of $\mathbb M_{d_n}(\mathbb C)$, for all $n\in\mathbb N$, such that $P\subset\prod_\omega P_n$ and $Q\subset\prod_\omega Q_n$.
\end{mcor}

By \cite[Theorem B]{IS19}, the conclusion of Corollary \ref{C} holds for the ultrapower $M^\omega$ of certain, fairly complicated, examples of II$_1$ factors $M$.
Corollary \ref{C} provides the first natural examples of tracial ultraproducts that satisfy its conclusion.
We conjecture that this phenomenon holds for any ultraproduct II$_1$ factor
 $\prod_\omega M_n$. 

Second, Theorem \ref{main} can be reformulated as a property of the full group C$^*$-algebra C$^*(\mathbb F_2\times\mathbb F_2)$.
This has been an important object of study since the work of Kirchberg \cite{Ki93} showing that certain properties of C$^*(\mathbb F_2\times\mathbb F_2)$ (being residually finite or having a faithful trace) are equivalent to  Connes' embedding problem (see \cite{Oz04,Pi20}).

Theorem \ref{main} implies the existence of a $*$-homomorphism $\varphi:\text{C}^*(\mathbb F_2\times\mathbb F_2)\rightarrow\prod_\omega \mathbb M_{d_n}(\mathbb C)$, for a sequence $(d_n)\subset\mathbb N$, which does not ``lift" to a $*$-homomorphism $\widetilde\varphi:\text{C}^*(\mathbb F_2\times\mathbb F_2)\rightarrow\ell^{\infty}(\mathbb N,\mathbb M_{d_n}(\mathbb C))$. Specifically, there is no $*$-homomorphism $\widetilde\varphi$ such that $\pi\circ\widetilde\varphi=\varphi$, where $\pi:\ell^{\infty}(\mathbb N,\mathbb M_{d_n})\rightarrow \prod_\omega \mathbb M_{d_n}$ is the quotient homomorphism.
We do not know if $\varphi$ admits a unital completely positive (ucp) lift $\widetilde\varphi$. If no ucp lift exists, then it would follow that $\text{C}^*(\mathbb F_2\times\mathbb F_2)$ does not have the local lifting property (LLP) (see \cite[Corollary 3.12]{Oz04}). Whether $\text{C}^*(\mathbb F_2\times\mathbb F_2)$ has the LLP is an open problem which goes back to \cite{Oz04} (see also \cite{Oz13,Pi20}).

\subsection*{Comments on the proof of Theorem \ref{main}} We end the introduction with a detailed outline of the proof of Theorem \ref{main}. Let us first reduce it to a simpler statement.
As we prove in Lemma \ref{flexible}, if $\Gamma_1,\Gamma_2$ are HS-stable, then $\Gamma_1\times\Gamma_2$ is flexibly HS-stable if and only if it is HS-stable. Also, if  $(\Gamma_1*\Lambda_1)\times(\Gamma_2*\Lambda_2)$ is HS-stable,  for groups $\Gamma_1,\Gamma_2,\Lambda_1,\Lambda_2$, then $\Gamma_1\times\Gamma_2$  must be HS-stable. These facts imply that proving Theorem \ref{main} is equivalent to showing that $\mathbb F_2\times\mathbb F_2$ is not HS-stable. To prove the latter statement, we will reason by contradiction assuming that $\mathbb F_2\times\mathbb F_2$ {\it is} HS-stable.

The proof of Theorem \ref{main} is divided into two parts, which we discuss separately below.
A main novelty of our approach is the use of ideas and techniques from the theory of (infinite dimensional) von Neumann algebras to prove a statement concerning finite unitary matrices. We combine small perturbations results for von Neumann algebras with finite dimensional analogues of two key ideas (the use of deformations and spectral gap arguments) from Popa's deformation/rigidity theory. 

{\it The first part of the proof}, which occupies Sections 3-5, is devoted to proving the following: 

\begin{proposition}\label{reduction}
Assume that $\mathbb F_2\times\mathbb F_2$ is \emph{HS}-stable. 
Then for every $\varepsilon>0$, there exists $\delta>0$ such that the following holds: 
for every $k,m,n\in\mathbb N$ and every $U_1,...,U_k,V_1,...,V_m\in \emph{U}(n)$ satisfying that $\frac{1}{km}\sum_{i=1}^k\sum_{j=1}^m\|[U_i,V_j]\|_2^2\leq\delta$, we can find  $\widetilde U_1,...,\widetilde U_k,\widetilde V_1,...,\widetilde V_m\in\emph{U}(n)$  such that 
\begin{enumerate}
\item $[\widetilde U_i,\widetilde V_j]=0$, for every $1\leq i\leq k$ and  $1\leq j\leq m$, 
\item $\frac{1}{k}\sum_{i=1}^k\|U_i-\widetilde U_i\|_2^2\leq\varepsilon$ and $\frac{1}{m}\sum_{j=1}^m\|V_j-\widetilde V_j\|_2^2\leq\varepsilon$.
\end{enumerate}
\end{proposition}
To illustrate the strength of the conclusion of Proposition \ref{reduction}, we make the following remark:
\begin{remark}\label{hs}
Let $k,m\in\mathbb N$. Then $\mathbb F_k\times\mathbb F_m$ is HS-stable if and only if for every  $\varepsilon>0$, there exists $\delta>0$ such that the following holds: for every $n\in\mathbb N$ and every $U_1,...,U_k,V_1,...,V_m\in \text{U}(n)$ satisfying that $\|[U_i,V_j]\|_2\leq\delta$, for all $1\leq i\leq k$ and $1\leq j\leq m$, we can find $\widetilde U_1,...,\widetilde U_k,\widetilde V_1,...,\widetilde V_m\in \text{U}(n)$  such that  $[\widetilde U_i,\widetilde V_j]=0$,   $\|U_i-\widetilde U_i\|_2\leq\varepsilon$ and  $\|V_j-\widetilde V_j\|_2\leq\varepsilon$, for all $1\leq i\leq k$ and $1\leq j\leq m$.
In view of this, Proposition \ref{reduction} can be interpreted as follows: if $\mathbb F_2\times\mathbb F_2$ is HS-stable then $\mathbb F_k\times\mathbb F_m$ is HS-stable and, moreover, it satisfies an ``averaged" version of HS-stability, {\it uniformly} over all $k,m\in\mathbb N$. 
\end{remark}

We continue with some comments on the proof of Proposition \ref{reduction} under the stronger assumption that $\mathbb F_3\times\mathbb F_3$ is HS-stable. 
The proof of Proposition \ref{reduction} has three main ingredients. All subalgebras of matrix algebras considered below are taken to be von Neumann (i.e., self-adjoint) subalgebras.

The first is a small perturbation result for subalgebras of a tensor product of three matrix algebras $M=\mathbb M_k(\mathbb C)\otimes\mathbb M_n(\mathbb C)\otimes\mathbb M_m(\mathbb C)$ (see Lemma \ref{perturb}). Roughly speaking, we prove that 
 any subalgebra $P\subset M$ which almost contains $\mathbb M_k(\mathbb C)\otimes 1\otimes 1$ and is almost contained in $\mathbb M_k(\mathbb C)\otimes\mathbb M_n(\mathbb C)\otimes 1$ must be close to a subalgebra of the form $\mathbb M_k(\mathbb C)\otimes S\otimes 1$, for some subalgebra $S\subset\mathbb M_n(\mathbb C)$.
Here, for subalgebras $P,Q\subset M$ and $\varepsilon>0$, we say that $P$ is $\varepsilon$-contained in $Q$ if for all $x\in P$ with $\|x\|\leq 1$ there is $y\in Q$ with $\|x-y\|_2\leq \varepsilon$, and that $P$ is $\varepsilon$-close to $Q$ if $P\subset_\varepsilon Q$ and $Q\subset_\varepsilon P$ \cite{MvN43}.  A crucial aspect of Lemma \ref{perturb} is that the constants involved are independent of $k,n,m\in\mathbb N$. Its proof is based on ideas from \cite{Ch79,Po01,Po03} and in particular uses the basic construction as in \cite{Ch79}.

The second ingredient in the proof of the Proposition \ref{reduction} is the existence of pairs of unitaries satisfying the following ``spectral gap" condition: for a universal constant $\kappa>0$ and every $n\in\mathbb N$, we can find $X_1,X_2\in\text{U}(n)$ such that $\|x-\tau(x)1\|_2\leq \kappa(\|[X_1,x]\|_2+\|[X_2,x]\|_2)$, for every $x\in\mathbb M_n(\mathbb C)$ (see Lemma \ref{gap}).
This is a consequence of a result of Hastings \cite{Ha07,Pi12} on quantum expanders. 

To finish the proof of Proposition \ref{reduction} we combine the first two ingredients with a tensor product trick. Let $U_1,\cdots,U_k,V_1,\cdots,V_m\in\text{U}(n)$ such that $\|[U_i,V_j]\|_2\approx 0$, for every $i,j$. Let $X_1,X_2\in \text{U}(k)$ and $Y_1,Y_2\in\text{U}(m)$ be pairs of unitaries with spectral gap.
We define $M=\mathbb M_k(\mathbb C)\otimes\mathbb M_n(\mathbb C)\otimes\mathbb M_m(\mathbb C)$ and unitaries $Z_1,Z_2,Z_3,T_1,T_2,T_3\in M$ by letting $$\text{$Z_1=X_1\otimes 1\otimes 1, \;\;Z_2=X_2\otimes 1\otimes 1, \;\; Z_3=\sum_{i=1}^ke_{i,i}\otimes U_i\otimes 1$}$$ $$ T_1=1\otimes 1\otimes Y_1,\;\; T_2=1\otimes 1\otimes Y_2,\;\; T_3=\sum_{j=1}^m1\otimes V_j\otimes e_{j,j}.$$
Then $\|[Z_i,T_j]\|_2\approx 0$, for every $1\leq i,j\leq 3$.  Since $\mathbb F_3\times\mathbb F_3$ is assumed HS-stable, there are unitaries $Z_i',T_j'\in M$ such that $[Z_i',T_j']=0$, $\|Z_i-Z_i'\|_2\approx 0$ and $\|T_j-T_j'\|_2\approx 0$, for every $1\leq i,j\leq 3$. 

Let $P$ the subalgebra of $M$ generated by $Z_1',Z_2',Z_3'$ and let $Q$ be its commutant.
Since $P$ almost commutes with $T_1,T_2$, using the spectral gap property of $Y_1,Y_2$ via an argument inspired by \cite{Po06a,Po06b} we deduce that $P$ is almost contained in $\mathbb M_{k}(\mathbb C)\otimes\mathbb M_n(\mathbb C)\otimes 1$. 
Similarly, it follows that $Q$ is almost contained in $1\otimes\mathbb M_n(\mathbb C)\otimes\mathbb M_m(\mathbb C)$. Since $P$ is the commutant of $Q$, the bicommutant theorem  implies that $\mathbb M_k(\mathbb C)\otimes 1\otimes 1$ is almost contained in $P$. The first ingredient of the proof now provides commuting subalgebras $R,S\subset \mathbb M_n(\mathbb C)$ such that $P$ is close to $\mathbb M_k(\mathbb C)\otimes R\otimes 1$ and $Q$ is close to $1\otimes S\otimes\mathbb M_m(\mathbb C)$. At this point, the conclusion of Proposition \ref{reduction} follows easily.

 In {\it the second part of the proof} of Theorem \ref{main}, presented in Section 6, we construct a counterexample to the conclusion of Proposition \ref{reduction} and derive that $\mathbb F_2\times\mathbb F_2$ is not HS-stable. 
 Our construction, which we describe in detail below, is inspired by Popa's malleable deformation for noncommutative Bernoulli actions, see \cite{Po03,Va06}, and its variant introduced in \cite{Io06}.
 \begin{construction}
 Let $n\in\mathbb N$ and $t\in\mathbb R$. 
\begin{enumerate}
\item We denote $M_n=\bigotimes_{k=1}^n\mathbb M_2(\mathbb C)\cong\mathbb M_{2^n}(\mathbb C)$ and $A_n=\bigotimes_{k=1}^n\mathbb C^2\cong \mathbb C^{2^n}$.  We view $A_n$ as a subalgebra of $M_n$, where we embed $\mathbb C^2\subset\mathbb M_2(\mathbb C)$ as the diagonal matrices.

 \item For $1\leq i\leq n$, let $X_{n,i}=1\otimes\cdots 1\otimes \sigma\otimes 1 \cdots\otimes 1\in A_n$, where $\sigma=\begin{pmatrix} 1&0\\ 0&-1\end{pmatrix}\in\mathbb C^2$ is placed on the $i$-th tensor position.

 \item Let $G_n$ be a finite group of unitaries which generates $A_n\otimes M_n$. 


 \item We define $U_t\in\text{U}(\mathbb C^2\otimes\mathbb C^2)$ by  $U_t=P+e^{it}(1-P)$, where $P:\mathbb C^2\otimes\mathbb C^2\rightarrow\mathbb C^2\otimes\mathbb C^2$ be the orthogonal projection onto the one dimensional space spanned by $e_1\otimes e_2-e_2\otimes e_1$ .

 \item We identify $M_n\otimes M_n=\otimes_{k=1}^n(\mathbb M_2(\mathbb C)\otimes\mathbb M_2(\mathbb C))$ and  let $\theta_{t,n}$ be the automorphism of $M_n\otimes M_n$ given by $\theta_{t,n}(\otimes_{k=1}^nx_k)=\otimes_{k=1}^nU_tx_kU_t^*$. 
\item Finally, consider the following two sets of unitaries in $M_n\otimes M_n$: $\mathcal U_n=\{X_{n,i}\otimes 1\mid 1\leq i\leq n\}$ and $\mathcal V_{t,n}=G_n\cup\theta_{t,n}(G_n)$.
\end{enumerate}
 \end{construction}

 Then $\mathcal U_n$ and $\mathcal V_{t,n}$ almost commute: $\|[U,V]\|_2\leq 4t$, for $U\in\mathcal U_n$, $V\in\mathcal V_{t,n}$. This is because $\mathcal U_n$ commutes with $G_n$ and $\|\theta_{t,n}(U)-U\|_2\leq 2t$, for every $U\in\mathcal U_n$.
 Using this, we show that if $t>0$ is small enough, then the sets $\mathcal U_n,\mathcal V_{t,n}$ contradict the conclusion of Proposition \ref{reduction} for large $n\in\mathbb N$.

To informally outline our argument, suppose that $\mathbb F_2\times\mathbb F_2$ is HS-stable. Then Proposition \ref{reduction} provides commuting subalgebras $P_n,Q_n$ of $M_n\otimes M_n$ so that $P_n$ almost contains $\mathcal U_n$ and $Q_n$ almost contains $\mathcal V_{t,n}$.  Thus, $P_n$ almost commutes with $\mathcal V_{t,n}$ and hence with the generating groups $G_n$, $\theta_{t,n}(G_n)$ of $A_n\otimes M_n$, $\theta_{t,n}(A_n\otimes M_n)$. By passing to commutants, we derive that $P_n$ is almost contained in both $A_n\otimes 1$ and $\theta_{t,n}(A_n\otimes 1)$. 
By perturbing $P_n$ slightly, we can in fact assume that $P_n$ is a subalgebra of $A_n\otimes 1$ which is almost contained in $\theta_{t,n}(A_n\otimes 1)$.

Assume for a moment that $n=\infty$ in the above construction. Then $(\theta_{t,\infty})_{t\in\mathbb R}$ recovers the malleable deformation of the Bernoulli action on the hyperfinite II$_1$ factor $M_{\infty}=\overline{\otimes}_{k=1}^{\infty}\mathbb M_2(\mathbb C)$, see \cite{Po03,Va06}. In this case, if a subalgebra $P$ of $M_\infty\otimes 1$ is almost contained in $\theta_{t,\infty}(M_\infty\otimes 1)$, then it must have a finite dimensional direct summand \cite{Io06}.  While this result cannot be used in our finite dimensional setting, we use the intuition behind its proof  and a dimension argument to derive a contradiction.

Since $P_n\subset A_n\otimes 1$ is almost contained in $\theta_{t,n}(A_n\otimes 1)$, it is almost contained in the subspace of tensors of $A_n=\otimes_{k=1}^n\mathbb C^2$ of length at most $l$, for some $l$ independent on $n$. This forces the dimension of $P_n$ to be at most polynomial in $n$. But $P_n$ also almost contains $\mathcal U_n$ and so all tensors of length $1$. This forces the dimension of $P_n$ to be at least exponential in $n$, giving a contradiction  as $n\rightarrow\infty$.


\begin{remark} Let $\Gamma=\mathbb F_k\times\mathbb F_m$, for integers $k,m\geq 2$. Note that  $\mathcal U_n$ and $\mathcal V_{t,n}$ almost commute in the operator norm: $\|[U,V]\|\leq 4t$, for $U\in\mathcal U_n, V\in\mathcal V_{t,n}$. Using this fact, a close inspection of the proof of Theorem \ref{main} shows that we prove the following stronger statement: the asymptotic homomorphism $\varphi_n:\Gamma\rightarrow\text{U}(d_n)$ which witnesses that $\Gamma$ is not flexibly HS-stable is an asymptotic homomorphism in the operator norm, i.e., $\lim\limits_{n\rightarrow\infty}\|\varphi_n(gh)-\varphi_n(g)\varphi_n(h)\|=0$, for all $g,h\in\Gamma$.  Therefore, $\Gamma$ fails a hybrid notion of stability which weakens both the notion of  matricial stability studied in  \cite{ESS18,Da20} and (flexible) HS-stability.

 
\end{remark}

\subsection{Acknowledgements} I am grateful to R\'{e}mi Boutonnet for stimulating discussions and to Sorin Popa for helpful comments.

\section{Preliminaries}

While the main results of this paper concern matrix algebras, 
the proofs are based on ideas and techniques from the theory of von Neumann algebras. 
 Moreover, our proofs often extend with no additional effort from matrix algebras to general tracial von Neumann algebras.  
As such, it will be convenient to work in the latter framework.
In this section, we recall several basic notions and constructions concerning  von Neumann algebras (see \cite{AP} and \cite{Ta79} for more information). 

\subsection{Von Neumann algebras}\label{prelim}
For a complex Hilbert space $\mathcal H$, we denote by $\mathbb B(\mathcal H)$ the algebra of bounded linear operators on $\mathcal H$ and by $\text{U}(\mathcal H)=\{u\in\mathbb B(\mathcal H)\mid u^*u=uu^*=1\}$ the group of unitary operators on $\mathcal H$. 
 For $x\in\mathbb B(\mathcal H)$, we denote by $\|x\|$ its operator  norm. 
  A set  of operators $S\subset\mathbb B(\mathcal H)$ is called {\it self-adjoint} if $x^*\in S$, for all $x\in S$.
 We denote by $S'$ the {\it commutant} of $S$, i.e., the set of operators $y\in\mathbb B(\mathcal H)$ such that $xy=yx$, for all $x\in S$. 

A self-adjoint subalgebra $M\subset\mathbb B(\mathcal H)$ is a {\it von Neumann algebra} if it is closed in the weak operator topology. By von Neumann bicommutant's theorem, a unital self-adjoint subalgebra $M\subset\mathbb B(\mathcal H)$  is a von Neumann algebra if and only if it is equal to its bicommutant, $M=(M')'$. 
From now on, we assume that all von Neumann algebras $M$ are unital.
We denote by $\mathcal Z(M)=M'\cap M$ the {\it center} of $M$, by  $(M)_1=\{x\in M\mid \|x\|\leq 1\}$ the {\it unit ball} of $M$, by $M_{+}=\{x\in M\mid x\geq 0\}$ the set of {\it positive} elements of $M$, and by $\mathcal U(M)$ the group of unitary operators in $M$. 
We call $M$ a {\it factor} if $\mathcal Z(M)=\mathbb C1$. 
Two projections $p,q\in M$ are Murray-von Neumann {\it equivalent} if there is a partial isometry $v\in M$ such that $v^*v=p$ and $vv^*=q$. A linear functional $\varphi:M\rightarrow\mathbb C$ is a called (a) a {\it state} if $\varphi(1)=1$ and $\varphi(x)\geq 0$, for every $x\in M_+$, (b) {\it faithful} if having $\varphi(x)=0$, for some $x\in M_{+}$, implies that $x=0$, and (c)
 {\it normal} if $\sup\varphi(x_i)=\varphi(\sup x_i)$, for any increasing net $(x_i)\subset M_{+}$.

A  {\it tracial von Neumann algebra} is a pair $(M,\tau)$ consisting of a von Neumann algebra $M$ and a {\it trace} $\tau$, i.e., a faithful normal state $\tau:M\rightarrow\mathbb C$ which satisfies $\tau(xy)=\tau(yx)$, for all $x,y\in M$. 
We endow  $M$ with the $1$- and $2$-norms given by $\|x\|_1=\tau((x^*x)^{\frac{1}{2}})$ and  $\|x\|_{2}=\tau(x^*x)^{\frac{1}{2}}$, 
for all $x\in M$. Then $\|xy\|_{2}\leq \|x\|\;\|y\|_{2}$ and $\|xy\|_{2}\leq \|x\|_{2}\;\|y\|$, for all $x,y\in M$. 
We denote by L$^2(M)$ the Hilbert space obtained as the closure of $M$ with respect to $\|\cdot\|_{2}$, and consider the standard representation $M\subset\mathbb B(\text{L}^2(M))$ given by the left multiplication action of $M$ on $\text{L}^2(M)$. 
For further reference, we recall the Powers-St{\o}rmer inequality (see \cite[Proposition 6.2.4]{BO08} and \cite[Theorem 7.3.7]{AP})
\begin{equation}\label{PS} \text{$\|h-k\|_2^2\leq \|h^2-k^2\|_1\leq\|h-k\|_2\|h+k\|_2$, \;\; for every  $h,k\in M_+$,} \end{equation} and the following inequality
\begin{equation}\label{pq}\text{$ |\tau(p)-\tau(q)|\leq \|p-q\|_{2}^2$, \;\; for every projections $p,q\in M$.}
\end{equation}
The latter inequality holds because $\|p-q\|_{2}^2=\tau(p)+\tau(q)-2\tau(pq)$ and $\tau(pq)\leq\min\{\tau(p),\tau(q)\}$.
We also note that \eqref{PS} and \eqref{pq} more generally hold when $\tau:M\rightarrow\mathbb C$ is a semifinite trace.

The matrix algebra $\mathbb M_n(\mathbb C)=\mathbb B(\mathbb C^n)$, for $n\in\mathbb N$, 
 with its {\it normalized trace} $\tau:\mathbb M_n(\mathbb C)\rightarrow\mathbb C$ given by $$\text{$\tau(x)=\frac{1}{n}\sum_{i=1}^nx_{i,i}$\;\; for\;\; $x=(x_{i,j})_{i,j=1}^n\in \mathbb M_n(\mathbb C),$}$$ is a tracial von Neumann algebra.
The associated $2$-norm is the
 {\it normalized Hilbert-Schmidt norm}  $$\text{$\|x\|_2= \Big(\frac{1}{n}\sum_{i,j=1}^n|x_{i,j}|^2\Big)^{\frac{1}{2}}$,\;\; for \;\; $x=(x_{i,j})_{i,j=1}^n\in \mathbb M_n(\mathbb C)$.}$$

Since $\mathbb M_n(\mathbb C)$ has trivial center, it is a tracial factor. Any tracial factor  
is either finite dimensional and isomorphic to  $\mathbb M_n(\mathbb C)$, for  some $n\in\mathbb N$, or infinite dimensional and called a {\it {\rm II}$_1$ factor}. 

Moreover, any finite dimensional von Neumann algebra $M$ is isomorphic to a direct sum of matrix algebras and therefore it is tracial. Indeed, if $z_1,\cdots,z_k\in\mathcal Z(M)$ are the minimal projections, then $M=\bigoplus_{i=1}^kMz_i$, where  $Mz_i$ is a finite dimensional factor and thus a matrix algebra, for $1\leq i\leq k$.
We claim that there is a finite subgroup $G\subset\mathcal U(M)$ which generates $M$. If $M=\mathbb M_n(\mathbb C)$, we can take $G$ to be group of unitaries of the form $(\varepsilon_i\delta_{\sigma(i),j})_{i,j=1}^n$, where $\varepsilon_1,\cdots,\varepsilon_n\in\{\pm 1\}$ and $\sigma$ is a permutation of $\{1,\cdots,n\}$. In general, if $G_i\subset\mathcal U(Mz_i)$ is a generating group, for every $1\leq i\leq k$, then the finite group $G=\{u_1\oplus\cdots\oplus u_k\mid u_1\in G_1,\cdots,u_k\in G_k\}$ generates $M$.

A subalgebra of a matrix algebra $\mathbb M_n(\mathbb C)$ is a von Neumann subalgebra if and only if it is self-adjoint.
Nevertheless, for consistency, we will call self-adjoint subalgebras of $\mathbb M_n(\mathbb C)$ von Neumann algebras.

\subsection{The basic construction}\label{basicc}
Let $(M,\tau)$ be a tracial von Neumann algebra 
together with a von Neumann subalgebra $Q\subset M$. 
Then we have an embedding $\text{L}^2(Q)\subset \text{L}^2(M)$. We denote by $e_Q:\text{L}^2(M)\rightarrow \text{L}^2(Q)$ the orthogonal projection onto $\text{L}^2(Q)$. 
We denote by $\text{E}_Q:M\rightarrow Q$ the {\it conditional expectation} onto $Q$, i.e., the unique map satisfying  $\tau(\text{E}_Q(x)y)=\tau(xy)$, for all $x\in M$ and $y\in Q$. Considering the natural embedding $M\subset\text{L}^2(M)$, then $e_Q(M)\subset Q$ and $\text{E}_Q={e_Q}_{|M}$.

 {\it Jones' basic construction} $\langle M,e_Q\rangle$ of the inclusion $Q\subset M$ is defined as the von Neumann subalgebra of $\mathbb B(\text{L}^2(M))$ generated by $M$ and $e_Q$. Let $J:\text{L}^2(M)\rightarrow\text{L}^2(M)$ be the involution given by $J(x)=x^*$, for all $x\in M$. Then $\langle M,e_Q\rangle$ is equal to both $JQ'J$, the commutant of the right multiplication action of $Q$ on $\text{L}^2(M)$, and the weak operator closure of the span of $\{xe_Qy\mid x,y\in M\}$.
 
 The basic construction admits a normal semifinite trace $\text{Tr}:\langle M,e_Q\rangle_{+}\rightarrow [0,+\infty]$ 
which satisfies $$\text{$\text{Tr}(xe_Qy)=\tau(xy)$, for every $x,y\in M$.}$$
It also admits a
  normal semifinite center-valued tracial weight $\Phi:\langle M,e_Q\rangle_{+}\rightarrow\widehat{\mathcal Z(Q)}_+$ which satisfies $\text{Tr}=\tau\circ\Phi$, $\Phi(X^*X)=\Phi(XX^*)$, for every $X\in\langle M,e_Q\rangle$, and $\Phi(ST^*)=\text{E}_{\mathcal Z(Q)}(T^*S)$, for all bounded right $Q$-linear operators $S,T:\text{L}^2(Q)\rightarrow\text{L}^2(M)$    (see \cite[Section 9.4]{AP}). 
   Here, $\widehat{\mathcal Z(Q)}_+$ denotes the  set of positive operators affiliated with $\mathcal Z(Q)$. 
  If $S$ and $T$ are the left multiplication operators by $x$ and $y^*$, for $x,y\in M$,  then $ST^*=xe_Qy$ and $T^*S=\text{E}_Q(yx)$. Thus, we conclude that \begin{equation}
 \label{Phi} \text{$\Phi(xe_Qy)=\text{E}_{\mathcal Z(Q)}(yx)$, for every $x,y\in M$.}
  \end{equation}
  If $z\in\mathcal Z(Q)$, then $e_QJzJ=e_Qz^*$ and thus $\Phi(xe_QyJzJ)=\Phi(xe_Qz^*y)=\text{E}_{\mathcal Z(Q)}(z^*yx)=\Phi(xe_Qy)z^*$, for every $x,y\in M$. Therefore, it follows that 
  \begin{equation}\label{bimod}
  \text{$\Phi(TJzJ)=\Phi(T)z^*$, for every $T\in\langle M,e_Q\rangle_+$ and $z\in \mathcal Z(Q)$.}
  \end{equation}

We note that if $M$ is finite dimensional, then  $\langle M,e_Q\rangle$ is finite dimensional and $\text{Tr}$ and $\Phi$ are finite, i.e., $\text{Tr}(1)<\infty$ and $\Phi(1)\in\mathcal Z(Q)$. 
We next record two well-known properties of $\Phi$.
\begin{lemma}\label{twoproj}
  Given two projections $p,q\in \langle M,e_Q\rangle$, the following hold:
  \begin{enumerate}
 \item $p$ is equivalent to a subprojection of $q$  if and only if $\Phi(p)\leq \Phi(q)$.
 \item There is a projection $r\in\langle M,e_Q\rangle$ such that $r\leq p$ and $\Phi(r)=\min\{\Phi(p),\Phi(q)\}$.
  \end{enumerate}
\end{lemma}

{\it Proof.} For (1), see \cite[Proposition 9.1.8.]{AP}. For (2), it is easy to see that any maximal projection $r\leq p$ which is equivalent to a subprojection of $q$ has the desired property.  
\hfill$\blacksquare$

\subsection{Almost containment} Let us recall the notion of $\varepsilon$-containment studied in \cite{MvN43,Mc70,Ch79}. Let  $P\subset pMp$, $Q\subset qMq$ be von Neumann subalgebras of a tracial von Neumann algebra $(M,\tau)$, for projections $p,q\in M$.
For $\varepsilon\geq 0$, we  write $P\subset_{\varepsilon}Q$ and say that $P$ is {\it $\varepsilon$-contained} in $Q$  if $\|x-\text{E}_Q(x)\|_{2}\leq\varepsilon$, for all $x\in (P)_1$. 
We also define the distance between $P$ and $Q$ by letting $${\bf d}(P,Q):=\min\{\varepsilon\geq 0\mid  P\subset_{\varepsilon}Q\;\;\;\text{and}\;\;\; Q\subset_{\varepsilon}P\}.$$

{\bf Convention.} To specify the trace $\tau$, we sometimes write $\|x\|_{2,\tau}$, $\subset_{\varepsilon,\tau}$, $\text{\bf d}_{\tau}$ instead of $\|x\|_2$, $\subset_{\varepsilon}$, $\text{\bf d}$.

In the rest of this subsection, we prove several useful lemmas.
We start with two well-known results.

\begin{lemma}\label{close} The following hold:
\begin{enumerate}
\item Let $M$ be a von Neumann algebra with a faithful normal semifinite trace $\tau$. If $p,q\in M$ are equivalent finite projections, there is $u\in \mathcal U(M)$ satisfying $upu^*=q$ and $\|u-1\|_2\leq 3\|p-q\|_2$.
\item Let $(M,\tau)$ be a tracial von Neumann algebra, $P\subset M$ a von Neumann subalgebra and $u\in \mathcal U(M)$. Then there is $v\in \mathcal U(P)$ satisfying $\|u-v\|_2\leq 3\|u-\emph{E}_P(u)\|_2$.
\end{enumerate}
\end{lemma}

{\it Proof.}
For (1), see \cite[Lemma 2.2]{Ch79}. For (2), using the polar decomposition of $\text{E}_P(u)$ we find $v\in\mathcal U(P)$ such that $\text{E}_P(u)=v|\text{E}_P(u)|$. 
Then $\|u-v\|_2\leq \|u-\text{E}_P(u)\|_2+\|1-|\text{E}_P(u)|\|_2$. Since $\|1-|\text{E}_P(u)|\|_2\leq \|1-|\text{E}_P(u)|^2\|_2=\|u^*u-\text{E}_P(u)^*\text{E}_P(u)\|_2\leq 2\|u-\text{E}_P(u)\|_2$, (2) follows.
\hfill$\blacksquare$

\begin{lemma}\label{commutant}
Let $(M,\tau)$ be a tracial von Neumann algebra and $P\subset pMp, Q\subset qMq$ be von Neumann subalgebras, for some projections $p,q\in M$. Assume that $P\subset_{\varepsilon}Q$ and $\|p-q\|_{2}\leq\varepsilon$, for some $\varepsilon>0$.
Then $Q'\cap qMq\subset_{4\varepsilon} P'\cap pMp$.
\end{lemma}

{\it Proof.}
Let $y\in (Q'\cap qMq)_1$. Then for every $u\in\mathcal U(P)$, we have that $[u,pyp]=p[u,y]p$ and thus
$\|[u,pyp]\|_{2}\leq \|[u,y]\|_{2}=\|[u-\text{E}_Q(u),y]\|_{2}\leq 2\|u-\text{E}_Q(u)\|_{2}\leq 2\varepsilon.$
Since this holds for every $u\in\mathcal U(P)$, it follows that $\|pyp-\text{E}_{P'\cap pMp}(pyp)\|_{2}\leq 2\varepsilon$. 
Since $y=qyq$, we also have that $\|y-pyp\|_{2}=\|qyq-pyp\|_{2}\leq 2\|p-q\|_{2}\leq 2\varepsilon$.  By combining the last two inequalities we derive that $\|y-\text{E}_{P'\cap pMp}(y)\|_{2}\leq\|y-\text{E}_{P'\cap pMp}(pyp)\|_{2}\leq 4\varepsilon$, which finishes the proof.
\hfill$\blacksquare$

The following lemma is a simple application of the basic construction.

\begin{lemma}\label{basic}
Let $(M,\tau)$ be a tracial von Neumann algebra and let $P,Q\subset M$ be von Neumann subalgebras. Assume that $P$ is finite dimensional and $G\subset\mathcal U(P)$ is a finite subgroup which generates $P$ and satisfies that $\frac{1}{|G|}\sum_{U\in G}\|U-\emph{E}_Q(U)\|_2^2\leq\varepsilon$, for some $\varepsilon>0$.
Then $P\subset_{\sqrt{2\varepsilon}}Q$.
\end{lemma}
{\it Proof.} Consider the basic construction $\langle M,e_Q\rangle$ with its canonical semifinite trace $\text{Tr}:\langle M,e_Q\rangle\rightarrow\mathbb C$. Denote $f=\frac{1}{|G|}\sum_{U\in G}Ue_QU^*\in\langle M,e_Q\rangle$. Then an easy calculation shows that \begin{equation}\label{feQ}
\|f-e_Q\|_{2,\text{Tr}}^2=\frac{1}{|G|}\sum_{U\in G}\|U-\text{E}_Q(U)\|_2^2\leq\varepsilon.
\end{equation}
Since $G$ is a group, $f$ commutes with $G$ and thus with $P$. Therefore, if $x\in (P)_1$, then using that $x$ commutes with $f$ and \eqref{feQ} we get that $$\|x-\text{E}_Q(x)\|_2^2=\frac{1}{2}\|xe_Q-e_Qx\|_{2,\text{Tr}}^2=\frac{1}{2}\|x(e_Q-f)-(e_Q-f)x\|_{2,\text{Tr}}^2\leq 2\|e_Q-f\|_{2,\text{Tr}}^2\leq 2\varepsilon.$$
This proves the conclusion.
\hfill$\blacksquare$

Our next goal is to establish the following useful elementary lemma.

\begin{lemma}\label{ccorner}
Let $(M,\tau)$ be a tracial factor and $P\subset M$ a von Neumann subalgebra. Let  $\varepsilon\in (0,\frac{1}{8}]$ and assume that $q\in M$ is a projection such that $\tau(1-q)=\varepsilon$. 
Then there is a von Neumann subalgebra $Q\subset qMq$ such that $\emph{\bf d}(P,Q)\leq 14\varepsilon^{\frac{1}{4}}.$
\end{lemma}
{\it Proof.} 
We claim that $P$ or $P'\cap M$ contains a projection $p$ such that $\tau(p)\in [\varepsilon,\varepsilon^{\frac{1}{2}}]$.
Assume that $P$ does not contain a projection $p$ with  $\tau(p)\in [\varepsilon,\varepsilon^{\frac{1}{2}}]$. 
Since $2\varepsilon<\varepsilon^{\frac{1}{2}}$, there is a minimal projection $r\in P$ such that $\tau(r)>\varepsilon^{\frac{1}{2}}$. 
  Let $z\in P$ be the smallest central projection such that $r\leq z$. Then we can find $d\in\mathbb N$ such that $\tau(z)=d\tau(r)$, $Pz\cong\mathbb M_d(\mathbb C)$ and there is a $*$-isomorphism $\theta:rMr\rightarrow z(P'\cap M)z$ such that $\tau(\theta(x))=d\tau(x)$, for every $x\in rMr$. Since $M$ is a factor and $\tau(1-q)=\varepsilon$,  $rMr$ contains a projection of trace $\varepsilon$. Thus, $P'\cap M$ contains a projection of trace $d\varepsilon$.
Since $d\varepsilon\geq \varepsilon$ and $d\varepsilon=(d\varepsilon^{\frac{1}{2}})\varepsilon^{\frac{1}{2}}\leq (d\tau(r))\varepsilon^{\frac{1}{2}}=\tau(z)\varepsilon^{\frac{1}{2}}\leq\varepsilon^{\frac{1}{2}}$, the claim follows.

Let $p$ be a projection in $P$ or $P'\cap M$ with $\tau(p)\in [\varepsilon,\varepsilon^{\frac{1}{2}}]$. Since $\tau(1-p)\leq 1-\varepsilon=\tau(q)$ and $M$ is a factor, there is projection $q_0\in qMq$ such that $\tau(q_0)=\tau(1-p)\geq 1-\varepsilon^{\frac{1}{2}}$. Then $1-p$ and $q_0$ are equivalent projections in $M$ such that $\|(1-p)-q_0\|_2\leq \|p\|_2+\|1-q_0\|_2=2\tau(p)^{\frac{1}{2}}\leq 2\varepsilon^{\frac{1}{4}}$. By Lemma \ref{close}(1) we can find a unitary $u\in M$ such that $u(1-p)u^*=q_0$ and $\|u-1\|_2\leq 6\varepsilon^{\frac{1}{4}}$.

Finally, since $1-p$ belongs to $P$ or $P'\cap M$, $(1-p)P(1-p)$ is a von Neumann algebra. Define $Q:=u(1-p)P(1-p)u^*\oplus\mathbb C(q-q_0).$ Then $Q\subset qMq$ is a von Neumann subalgebra and for
 $x\in (P)_1$ $$\|x-u(1-p)x(1-p)u^*\|_2\leq 2\|1-u(1-p)\|_2\leq 2(\|1-u\|_2+\|p\|_2)\leq 14\varepsilon^{\frac{1}{4}}.$$ Since 
$u(1-p)x(1-p)u^*\in Q$, we get that $P\subset_{14\varepsilon^{\frac{1}{4}}}Q$. Conversely, let $y\in (Q)_1$ and write $y=uxu^*+\alpha (q-q_0)$, for some $x\in (P)_1$ and $\alpha\in\mathbb C$ with $|\alpha|\leq1$. Then $$\|y-x\|_2\leq \|uxu^*-x\|_2+\|q-q_0\|_2\leq 2\|u-1\|_2+\|1-q_0\|_2\leq 13\varepsilon^{\frac{1}{4}}.$$ Hence, $Q\subset_{13\varepsilon^{\frac{1}{4}}}P$, and the conclusion follows.
\hfill$\blacksquare$

We end this section by illustrating the usefulness of Lemma \ref{ccorner} in proving the following:

\begin{lemma}\label{flexible} Let $\Gamma_1$ and $\Gamma_2$ be \emph{HS}-stable countable groups.
 Then $\Gamma_1\times\Gamma_2$ is \emph{HS}-stable if and only if it is flexibly \emph{HS}-stable.
\end{lemma}
{\it Proof.} Let $\Gamma=\Gamma_1\times\Gamma_2$. To prove the lemma, we only have to argue that if $\Gamma$ is flexibly HS-stable, then it is HS-stable. 
To this end, assume that $\Gamma$ is flexibly HS-stable.

For $n\in\mathbb N$, denote by $\tau_n$ the normalized trace of $\mathbb M_n(\mathbb C)$.
 Let $\pi_n:\Gamma\rightarrow \text{U}(k_n)$ be an asymptotic homomorphism. Since $\Gamma$ is flexibly HS-stable, we can find homomorphisms $\rho_n:\Gamma\rightarrow \text{U}(K_n)$, for $K_n\geq k_n$, such that $\lim_{n\rightarrow\infty}\frac{K_n}{k_n}=1$ and denoting by $e_n:\mathbb C^{K_n}\rightarrow\mathbb C^{k_n}$ the orthogonal projection, \begin{equation}\label{rhon}\text{$\lim\limits_{n\rightarrow\infty}\|\pi_n(g)-e_n\rho_n(g)e_n\|_{2,\tau_{k_n}}=0$, for every $g\in\Gamma$.}\end{equation} 

Since $\lim_{n\rightarrow\infty}\frac{K_n}{k_n}=1$ and $\|e_n-1\|_{2,\tau_{K_n}}=\sqrt{\frac{K_n-k_n}{K_n}}$, for every $n\in\mathbb N$, we get that \begin{equation}\label{e_n}\lim\limits_{n\rightarrow\infty}\|e_n-1\|_{2,\tau_{K_n}}=0.\end{equation}
Let $P_n\subset\mathbb M_{K_n}(\mathbb C)$ be the von Neumann algebra generated by $\rho_n(\Gamma_1)$. 
By using \eqref{e_n} and applying Lemma \ref{ccorner} we can find a von Neumann subalgebra $Q_n\subset e_n\mathbb M_{K_n}(\mathbb C)e_n\equiv \mathbb M_{k_n}(\mathbb C)$ such that \begin{equation}\label{PQ}\lim\limits_{n\rightarrow\infty}\text{\bf d}_{\tau_{K_n}}(P_n,Q_n)=0.\end{equation} Next, by combining \eqref{e_n}, \eqref{PQ} and  Lemma \ref{commutant} we derive that \begin{equation}\label{P'Q'}\lim\limits_{n\rightarrow\infty}\text{\bf d}_{\tau_{K_n}}(P_n'\cap\mathbb M_{K_n}(\mathbb C),Q_n'\cap \mathbb M_{k_n}(\mathbb C))=0.\end{equation} 
Since $\rho_n(\Gamma_1)\subset P_n$ and $\rho_n(\Gamma_2)\subset P_n'\cap\mathbb M_{K_n}(\mathbb C)$, combining \eqref{rhon}, \eqref{e_n}, \eqref{PQ} and \eqref{P'Q'} implies that 
$\lim_{n\rightarrow\infty}\|\pi_n(g_1)-\text{E}_{Q_n}(\pi_n(g_1))\|_{2,\tau_{k_n}}=0$ and
 $\lim_{n\rightarrow\infty}\|\pi_n(g_2)-\text{E}_{Q_n'\cap\mathbb M_{k_n}(\mathbb C)}(\pi_n(g_2))\|_{2,\tau_{k_n}}=0$, for every $g_1\in\Gamma_1$ and $g_2\in\Gamma_2$.
By putting together these facts and Lemma \ref{close}(2), we can find maps $\sigma_n^1:\Gamma_1\rightarrow\mathcal U(Q_n)$ and $\sigma_n^2:\Gamma_2\rightarrow\mathcal U(Q_n'\cap\mathbb M_{k_n}(\mathbb C))$ such that $\lim_{n\rightarrow\infty}\|\pi_n(g_1)-\sigma_n^1(g_1)\|_{2,\tau_{k_n}}=0$ and $\lim_{n\rightarrow\infty}\|\pi_n(g_2)-\sigma_n^2(g_2)\|_{2,\tau_{k_n}}=0$, for every $g_1\in\Gamma_1$ and $g_2\in\Gamma_2$.

Then $(\sigma_n^1)$ and $(\sigma_n^2)$ are asymptotic homomorphisms of $\Gamma_1$ and $\Gamma_2$, respectively.
Since $\Gamma_1$ and $\Gamma_2$ are HS-stable, we can find homomorphisms $\delta_n^1:\Gamma_1\rightarrow\mathcal U(Q_n)$ and $\delta_n^2:\Gamma_2\rightarrow\mathcal U(Q_n'\cap\mathbb M_{k_n}(\mathbb C))$ such that $\lim_{n\rightarrow\infty}\|\sigma_n^1(g_1)-\delta_n^1(g_1)\|_{2,\tau_{k_n}}=0$ and $\lim_{n\rightarrow\infty}\|\sigma_n^2(g_2)-\delta_n^2(g_2)\|_{2,\tau_{k_n}}=0$, for every $g_1\in\Gamma_1$ and $g_2\in\Gamma_2$.
Then we have $\lim_{n\rightarrow\infty}\|\pi_n(g_1)-\delta_n^1(g_1)\|_{2,\tau_{k_n}}=0$ and $\lim_{n\rightarrow\infty}\|\pi_n(g_2)-\delta_n^2(g_2)\|_{2,\tau_{k_n}}=0$, for every $g_1\in\Gamma_1$ and $g_2\in\Gamma_2$.

Finally, since $\delta_n^1$ and $\delta_n^2$ have commuting images for every $n\in\mathbb N$, we can define a homomorphism $\delta_n:\Gamma\rightarrow\text{U}(k_n)$ by letting $\delta_n(g_1,g_2)=\delta_n^1(g_1)\delta_n^2(g_2)$, for every $g_1\in\Gamma_1$ and $g_2\in\Gamma_2$. It is then clear that $\lim_{n\rightarrow\infty}\|\pi_n(g)-\delta_n(g)\|_{2,\tau_{k_n}}=0$, for every $g\in\Gamma$. This proves that $\Gamma$ is HS-stable.
\hfill$\blacksquare$

\section{Perturbation results} 
In this section we study the almost containment relation for tracial von Neumann algebras. A crucial feature of the results is that they do not depend on the dimensions of the algebras involved.
\subsection{A ``small perturbation" lemma}
Our main result is the following small perturbation lemma. If $P$ and $Q$ are subalgebras of a tracial von Neumann algebra such that $P$ is almost contained in $Q$, we show that $P$ must be close to a subalgebra of $\mathbb M_2(\mathbb C)\otimes Q$.
  For a tracial von Neumann algebra $(M,\tau)$, we equip $\mathbb M_2(\mathbb C)\otimes M$ with the trace $\widetilde\tau$ given by  $\widetilde\tau(\sum_{i,j=1}^2e_{i,j}\otimes x_{i,j})=\frac{1}{2}\big(\tau(x_{1,1})+\tau(x_{2,2})\big)$.

\begin{lemma}
\label{Ch}
Let $(M,\tau)$ be a tracial von Neumann algebra and  let $P,Q\subset M$ be von Neumann subalgebras. Assume that $P\subset_{\varepsilon}Q$, for some $\varepsilon\in (0,\frac{1}{200})$.  
Then there exists a $*$-homomorphism $\theta:P\rightarrow\mathbb M_2(\mathbb C)\otimes Q$ such that $\|\theta(x)-e_{1,1}\otimes x\|_{2,\widetilde\tau}\leq 30\varepsilon^{\frac{1}{8}}$, for every $x\in (P)_1$.
\end{lemma}

Lemma \ref{Ch} implies that there is a non-trivial $*$-homomorphism from $P$ to $\mathbb M_2(\mathbb C)\otimes Q$. This generalizes \cite[Theorem 4.7]{Ch79} where the same conclusion was proved assuming that $Q$ is a factor. 
As shown in \cite[Theorem A.2]{Po01} under certain conditions (e.g., if $P,Q\subset M$ are irreducible subfactors and $P\subset M$ is regular)  $P\subset_{\varepsilon}Q$ implies the existence of $u\in \mathcal U(M)$ such that $uPu^*\subset Q$.  However, such a strong conclusion does not hold in general even for irreducible subfactors $P,Q\subset M$ (see \cite[Proposition 5.5]{PSS03}).

The proof of Lemma \ref{Ch} relies on ideas of Christensen \cite{Ch79} and Popa \cite{Po01,Po03}. As in \cite{Ch79} we use the basic construction $\langle M,e_Q\rangle$ and find a projection $f\in P'\cap \langle M,e_Q\rangle$ which is close to $e_Q$. Then, inspired by an argument in \cite{Po03}, we show that after replacing $f$ by $fJzJ$, for a projection $z\in\mathcal Z(Q)$ close to $1$, one may assume that $e_{1,1}\otimes f$ is subequivalent to $1\otimes e_Q$ in  $\mathbb M_2(\mathbb C)\otimes\langle M,e_q\rangle$.

{\it Proof.} 
Let $\text{Tr}:\langle M,e_Q\rangle\rightarrow\mathbb C$  be the canonical tracial weight. Since
$$\|ue_Qu^*-e_Q\|_{2,\text{Tr}}^2=2(1-\text{Tr}(ue_Qu^*e_Q))=2(1-\tau(u\text{E}_Q(u)^*))=2\|u-\text{E}_Q(u)\|_2^2,$$ we get that
\begin{equation}\label{almost}\text{$\|ue_Qu^*-e_Q\|_{2,\text{Tr}}\leq \sqrt{2}\varepsilon$, for every $u\in\mathcal U(P)$.}
\end{equation}
Let $\mathcal C\subset\langle M,e_Q\rangle$ be the weak operator closure of the convex hull of the set $\{ue_Qu^*\mid u\in\mathcal U(P)\}$.  Then $\mathcal C$ is $\|\cdot\|_{2,\text{Tr}}$-closed and admits an element $h$ of minimal
$\|\cdot\|_{2,\text{Tr}}$-norm which satisfies $0\leq h\leq 1$, $h\in P'\cap \langle M,e_Q\rangle$ and $\|h-e_Q\|_{2,\text{Tr}}\leq\sqrt{2}\varepsilon$ by \eqref{almost} (see \cite[Section 3]{Ch79} or \cite[Lemma 14.3.3]{AP17}).

Define the spectral projection $f:={\bf 1}_{[1-\delta^{\frac{1}{2}},1]}(h)$, where $\delta=\sqrt{2}\varepsilon$.   Then $f\in P'\cap \langle M,e_Q\rangle$ and $\|f-e_Q\|_{2,\text{Tr}}\leq \delta^{\frac{1}{2}}(1-\delta^{\frac{1}{2}})^{-1}$ by \cite[Lemma 2.1]{Ch79}.
As $\varepsilon<\frac{1}{200}$, we have $(1-\delta^{\frac{1}{2}})^{-1}<\sqrt{2}$ and thus \begin{equation}\label{ef}\|f-e_Q\|_{2,\text{Tr}}\leq 2\varepsilon^{\frac{1}{2}}.\end{equation}
Since $\|e_Q\|_{2,\text{Tr}}=1$ and $2\varepsilon^{\frac{1}{2}}<1$, we get that $\|f+e_Q\|_{2,\text{Tr}}\leq 3$. Combining the Powers-St{\o}rmer inequality \eqref{PS} and \eqref{ef} 
we further get that \begin{equation}\label{ef1}\|f-e_Q\|_{1,\text{Tr}}\leq \|f-e_Q\|_{2,\text{Tr}}\cdot\|f+e_Q\|_{2,\text{Tr}}\leq 6\varepsilon^{\frac{1}{2}}.
\end{equation}

Let $\Phi:\langle M,e_Q\rangle_+\rightarrow\widehat{\mathcal Z(Q)}_+$ be the center-valued tracial weight defined in Subsection \ref{basicc}.
Define the spectral projection $z:=\text{\bf 1}_{[0,\varepsilon^{\frac{1}{4}}]}(\Phi(|f-e_Q|))$. Then $z\in\mathcal Z(Q)$ and $\varepsilon^{\frac{1}{4}}(1-z)\leq \Phi(|f-e_Q|)$. Since $\text{Tr}=\tau\circ\Phi$, \eqref{ef1} implies that  $\tau(\Phi(|f-e_Q|))=\|f-e_Q\|_{1,\text{Tr}}\leq 6\varepsilon^{\frac{1}{2}}$. Thus, $\tau(1-z)\leq 6\varepsilon^{\frac{1}{4}}$ and hence 
\begin{equation}\label{zz}\|1-z\|_{2,\tau}\leq 3\varepsilon^{\frac{1}{8}}.\end{equation}

Let $J:\text{L}^2(M)\rightarrow\text{L}^2(M)$ be the canonical involution given by $J(x)=x^*$ and put $z'=JzJ$. Since $\langle M,e_Q\rangle=JQ'J$, we have that 
 $z'\in M'\cap \mathcal Z(\langle M,e_Q\rangle)$. Thus, $g=fz'\in P'\cap \langle M,e_Q\rangle$ is a projection.
 Moreover, we have that $e_Qz'=e_Qz$, while \eqref{bimod} gives that $\Phi(x)z=\Phi(xz')$, for every $x\in \langle M,e_Q\rangle$. 
 Altogether, since $\Phi(|f-e_Q|)z\leq\varepsilon^{\frac{1}{4}}z$,  we get that
 $$\Phi(|g-e_Qz|)=\Phi(|fz'-e_Qz'|)=\Phi(|f-e_Q|z')=\Phi(|f-e_Q|)z\leq\varepsilon^{\frac{1}{4}}z=\varepsilon^{\frac{1}{4}}\Phi(e_Qz).$$
 
From this we deduce that \begin{equation}\label{bounds}(1-\varepsilon^{\frac{1}{4}})\Phi(e_Qz)\leq \Phi(g)\leq (1+\varepsilon^{\frac{1}{4}})\Phi(e_Qz).\end{equation}
By Lemma \ref{twoproj}(2) we can find a projection $p_1\in\langle M,e_Q\rangle$ such that $p_1\leq g$ and $\Phi(p_1)=\min\{\Phi(g),\Phi(e_Qz)\}$. Put $p_2=g-p_1$.
Since $\Phi(p_1)\leq\Phi(e_Qz)$, by Lemma \ref{twoproj}(1), there is a projection $q_1\in \langle M,e_Q\rangle$ such that $q_1\leq e_Qz$ and $q_1$ is equivalent to $p_1$. Since $\Phi(g-p_1)=\max\{0,\Phi(g)-\Phi(e_Qz)\}$, \eqref{bounds} implies that $\Phi(g-p_1)\leq\varepsilon^{\frac{1}{4}}\Phi(e_Qz)$ and thus $\text{Tr}(g-p_1)\leq \varepsilon^{\frac{1}{4}}$. Similarly, since $p_1$ and $q_1$ are equivalent we have $\Phi(p_1)=\Phi(q_1)$, thus $\Phi(e_Qz-q_1)=\Phi(e_Qz-p_1)\leq \varepsilon^{\frac{1}{4}}\Phi(e_Qz)$ and hence $\text{Tr}(e_Qz-q_1)\leq\varepsilon^{\frac{1}{4}}$. 

The last paragraph gives that $\|g-p_1\|_{2,\text{Tr}}\leq\varepsilon^{\frac{1}{8}}$ and $\|e_Qz-q_1\|_{2,\text{Tr}}\leq\varepsilon^{\frac{1}{8}}$. Moreover, \eqref{ef} gives that \begin{equation}\label{gg}\|g-e_Qz\|_{2,\text{Tr}}=\|(f-e_Q)z'\|_{2,\text{Tr}}\leq\|f-e_Q\|_{2,\text{Tr}}\leq 6\varepsilon^{\frac{1}{2}}.\end{equation} As $\varepsilon<\frac{1}{200}$, we have that $6\varepsilon^{\frac{1}{2}}<\varepsilon^{\frac{1}{8}}$ and the triangle inequality gives that $\|p_1-q_1\|_{2,\text{Tr}}\leq 3\varepsilon^{\frac{1}{8}}$. By \cite[Lemma 2.2]{Ch79}, there is a partial isometry $v_1\in\langle M,e_Q\rangle$ such that $v_1v_1^*=p_1,v_1^*v_1=q_1$ and $\|v_1-p_1\|_{2,\text{Tr}}\leq 6\|p_1-q_1\|_{2,\text{Tr}}\leq 18\varepsilon^{\frac{1}{8}}$. Since $\text{Tr}(p_2)=\text{Tr}(g-p_1)\leq \varepsilon^{\frac{1}{4}}$, we get $\|p_2\|_{2,\text{Tr}}\leq\varepsilon^{\frac{1}{8}}$ and so
\begin{equation}\label{v_1}
\|v_1-g\|_{2,\text{Tr}}\leq \|v_1-p_1\|_{2,\text{Tr}}+\|p_2\|_{2,\text{Tr}}\leq 19\varepsilon^{\frac{1}{8}}.
\end{equation}

Next, since $\Phi(p_2)=\Phi(g-p_1)\leq\varepsilon^{\frac{1}{4}}\Phi(e_Qz)\leq\Phi(e_Qz)$, by Lemma \ref{twoproj}(1)
we can find a partial isometry $v_2\in\langle M,e_Q\rangle$ with $v_2v_2^*=p_2$ and $v_2^*v_2\leq e_Q$. Then 
\begin{equation}\label{v_2}\|v_2\|_{2,\text{Tr}}=\|p_2\|_{2,\text{Tr}}\leq \varepsilon^{\frac{1}{8}}.\end{equation}

Let $v=\begin{pmatrix} v_1 & v_2\end{pmatrix}\in\mathbb M_{1,2}(\mathbb C)\otimes \langle M,e_Q\rangle$. Then $v$ is a partial isometry with $v^*v\in \mathbb M_2(\mathbb C)\otimes\langle M,e_Q\rangle$, $v^*v\leq 1\otimes e_Q$ and $vv^*=g$. 
As $g\in P'\cap\langle M,e_Q\rangle$, the map $P\ni x\mapsto v^*xv\in\mathbb M_2(\mathbb C)\otimes\langle M,e_Q\rangle$ is a $*$-homomorphism. If $x\in P$, then $v^*xv\in (1\otimes e_Q)(\mathbb M_2(\mathbb C)\otimes\langle M,e_Q\rangle)(1\otimes e_Q)=\mathbb M_2(\mathbb C)\otimes Qe_Q$. As $Q\ni y\mapsto ye\in Qe_Q$ is a $*$-isomorphism, there is a $*$-homomorphism $\theta:P\rightarrow\mathbb M_2(\mathbb C)\otimes Q$ such that \begin{equation}\label{theta} \text{$\theta(x)(1\otimes e_Q)=v^*xv$, for every $x\in P$.}
\end{equation}

 Since  $6\varepsilon^{\frac{1}{2}}<\varepsilon^{\frac{1}{8}}$, \eqref{zz} and \eqref{gg} give that $\|g-e_Q\|_{2,\text{Tr}}\leq \|g-e_Qz\|_{2,\text{Tr}}+\|1-z\|_{2,\tau}\leq 4\varepsilon^{\frac{1}{8}}$. 
Let $x\in (P)_1$. Since $g$ commutes with $x$ we have $gxg=xg$ and the above inequality and \eqref{v_1} imply
\begin{equation}\label{v1e_Q} \|v_1^*xv_1-xe_Q\|_{2,\text{Tr}}\leq\|v_1^*xv_1-gxg\|_{2,\text{Tr}}+\|xg-xe_Q\|_{2,\text{Tr}}\leq 2\|v_1-g\|_{2,\text{Tr}}+\|g-e_Q\|_{2,\text{Tr}}\leq 42\varepsilon^{\frac{1}{8}}.
\end{equation}

Finally, for $y=\sum_{i,j}^2e_{i,j}\otimes y_{i,j}\in\mathbb M_2(\mathbb C)\otimes \langle M,e_Q\rangle$, we denote $\widetilde{\text{Tr}}(y)=\frac{1}{2}(\text{Tr}(y_{1,1})+\text{Tr}(y_{2,2}))$ and $\|y\|_{2,\widetilde{\text{Tr}}}=(\widetilde{\text{Tr}}(y^*y))^{\frac{1}{2}}$.  Then $\|z(1\otimes e_Q)\|_{2,\widetilde{\text{Tr}}}=\|z\|_{2,\widetilde{\tau}}$, for every $z\in\mathbb M_2(\mathbb C)\otimes M$.
This fact together with \eqref{v_2}, \eqref{theta} and \eqref{v1e_Q} gives that for every $x\in (P)_1$ we have
\begin{align*}\|\theta(x)-e_{1,1}\otimes x\|_{2,\widetilde{\tau}}&=\|\theta(x)(1\otimes e_Q)-e_{1,1}\otimes xe_Q\|_{2,\widetilde{\text{Tr}}}\\&=
\|v^*xv-e_{1,1}\otimes xe_Q\|_{2,\widetilde{\text{Tr}}}\\&=\Big(\frac{\|v_1^*xv_1-xe_Q\|_{2,\text{Tr}}^2+\|v_1^*xv_2\|_{2,\text{Tr}}^2+\|v_2^*xv_1\|_{2,\text{Tr}}^2+\|v_2^*xv_2\|_{2,\text{Tr}}^2}{2}\Big)^{\frac{1}{2}}\\&
\leq 30\varepsilon^{\frac{1}{8}}.\end{align*} which finishes the proof.
\hfill$\blacksquare$

\subsection{From almost containment to containment} Let $P$ and $Q$ be von Neumann subalgebras of a tracial von Neumann algebra. If $P$ is  close to a subalgebra of $Q$, then $P$ is almost contained in $Q$. In this subsection, we use Lemma \ref{Ch} to prove that the converse holds provided that $Q$ is a factor (see Corollary \ref{Chr})  or a finite dimensional abelian algebra (see Corollary \ref{abelian}).

\begin{corollary}
\label{Chr} For any $\varepsilon>0$, there is $\delta=\delta_1(\varepsilon)>0$ such that the following holds. 
Let $(M,\tau)$ be a tracial von Neumann algebra and $P,Q\subset M$ be von Neumann subalgebras such that $P\subset_{\delta}Q$. Assume that $Q$ is a factor.  
Then there exists a von Neumann subalgebra $R\subset Q$ such that $\emph{\bf{d}}(P,R)\leq \varepsilon$.
\end{corollary}

{\it Proof.} Given $\varepsilon>0$, we will prove that any $\delta>0$ such that $\delta<10^{-16}$ and $400\delta^{\frac{1}{16}}<\varepsilon$ works. 
Assume that $P\subset_{\delta}Q$. 
Then Theorem \ref{Ch} gives a $*$-homomorphism $\theta:P\rightarrow\mathbb M_2(\mathbb C)\otimes Q$ with $\|\theta(x)-e_{1,1}\otimes x\|_{2,\widetilde\tau}\leq 30\delta^{\frac{1}{8}}$, for all $x\in (P)_1$.  Let $q=\theta(1)$ and $A=\theta(P)$.  Then $\|e_{1,1}\otimes 1-q\|_{2,\widetilde{\tau}}\leq 30\delta^{\frac{1}{8}}$ and $\emph{\bf d}_{\widetilde\tau}(e_{1,1}\otimes P,A)\leq 30\delta^{\frac{1}{8}}$.

Since $Q$ is a factor, we can find projections $r\in q(\mathbb M_2(\mathbb C)\otimes Q)q$ and $s\in e_{1,1}\otimes Q$ such that $\widetilde{\tau}(r)=\widetilde\tau(s)=\min\{\widetilde{\tau}(q),\frac{1}{2}\}$. 
By \eqref{pq}, we get that $\widetilde{\tau}(q-r)\leq |\widetilde{\tau}(q)-\frac{1}{2}|\leq \|q-e_{1,1}\otimes 1\|_{2,\widetilde\tau}^2\leq 900\delta^{\frac{1}{4}}$, so $\|q-r\|_{2,\widetilde\tau}\leq 30\delta^{\frac{1}{8}}$.   In particular, $\widetilde\tau(q)\geq\frac{1}{2}-900\delta^{\frac{1}{4}}\ge\frac{1}{3}$.
Similarly, $\|e_{1,1}\otimes 1-s\|_{2,\widetilde\tau}\leq 30\delta^{\frac{1}{8}}$, thus \begin{equation}\label{sr}\|s-r\|_{2,\widetilde\tau}\leq \|s-e_{1,1}\otimes 1\|_{2,\widetilde\tau}+\|e_{1,1}\otimes 1-q\|_{2,\widetilde\tau}+\|q-r\|_{2,\widetilde\tau}\leq 90\delta^{\frac{1}{8}}.\end{equation}

By Lemma \ref{ccorner}, we can find a von Neumann subalgebra $B\subset r(\mathbb M_2(\mathbb C)\otimes Q)r$ such that \begin{equation}\label{AB'}{\text{\bf{d}}}_{\widetilde\tau}(A,B)\leq 14\Big(\frac{\widetilde{\tau}(q-r)}{\widetilde\tau(q)}\Big)^{\frac{1}{4}}\leq 200\delta^{\frac{1}{16}}.\end{equation}

Since $\widetilde\tau(r)=\widetilde\tau(s)$ and $Q$ is a factor, using Lemma \ref{close}(1) and \eqref{sr} we find a unitary $u\in\mathbb M_2(\mathbb C)\otimes Q$ such that $uru^*=s$ and $\|u-1\|_{2,\widetilde\tau}\leq 3\|s-r\|_{2,\widetilde\tau}\leq 270\delta^{\frac{1}{8}}$. Let $C=uBu^*\oplus\mathbb C(e_{1,1}\otimes 1-s)$. Then $C\subset e_{1,1}\otimes Q$ is a von Neumann subalgebra such that  \begin{equation}\label{BC}{\text{\bf{d}}}_{\widetilde\tau}(B,C)\leq 2\|u-1\|_{2,\widetilde\tau}+\|e_{1,1}\otimes 1-s\|_{2,\widetilde\tau}\leq 570\delta^{\frac{1}{8}}.\end{equation}

If $R\subset Q$ is a von Neumann subalgebra such that $C=e_{1,1}\otimes R$, then \eqref{AB'} and \eqref{BC} imply that \begin{align*}\text{\bf{d}}_{\tau}(P,R)=\sqrt{2}\;\text{\bf{d}}_{\widetilde\tau}(e_{1,1}\otimes P,C)&\leq\sqrt{2}(\text{\bf{d}}_{\widetilde\tau}(e_{1,1}\otimes P,A)+\text{\bf{d}}_{\widetilde\tau}(A,B)+\text{\bf{d}}_{\widetilde\tau}(B,C))\\&\leq\sqrt{2}(30\delta^{\frac{1}{8}}+200\delta^{\frac{1}{16}}+570\delta^{\frac{1}{8}})\leq 400\delta^{\frac{1}{16}}.
\end{align*}
This finishes the proof of the lemma.
\hfill$\blacksquare$

\begin{corollary}\label{abelian}
For any $\varepsilon>0$, there is $\delta=\delta_2(\varepsilon)>0$ such that the following holds. 
Let $(M,\tau)$ be a tracial von Neumann algebra and $P,Q\subset M$ be finite dimensional von Neumann subalgebras such that $P\subset_{\delta}Q$. Assume that $Q$ is abelian.   
Then there exists a von Neumann subalgebra $R\subset Q$ such that $\emph{\bf{d}}(P,R)\leq\varepsilon$.
\end{corollary}

{\it Proof.} Given $\varepsilon>0$, we will prove that any $\delta>0$ such that $\delta<\frac{1}{200}$ and $200\delta^{\frac{1}{8}}<\varepsilon$ works. Assume that $P\subset_{\delta}Q$.
 We  will first show that $P$ has a large abelian direct summand.
 
Let $z\in\mathcal Z(P)$ be the largest projection such that $Pz$ is abelian. Since $P(1-z)$ has no abelian direct summand, we can find a projection $p\in P(1-z)$ with $\tau(p)\geq\frac{\tau(1-z)}{3}$  and a unitary $u\in P(1-z)$ such that  $p$ and $upu^*$ are orthogonal. Thus,  $\|[u,p]\|_{2,\tau}=\sqrt{2}\|p\|_{2,\tau}\geq \sqrt{\frac{2}{3}}\|1-z\|_{2,\tau}$. On the other hand, since $P\subset_{\delta}Q$ and $Q$ is abelian, we get that $\|[u,p]\|_{2,\tau}\leq 2(\|u-\text{E}_Q(u)\|_{2,\tau}+\|p-\text{E}_Q(p)\|_{2,\tau})\leq 4\delta$.  By combining the last two inequalities, we derive that \begin{equation}\label{zzz}\|1-z\|_{2,\tau}\leq 4\sqrt{\frac{3}{2}}\delta\leq 5\delta.\end{equation}

Let $\{p_i\}_{i=1}^m$ and $\{q_j\}_{j=1}^n$ be the minimal projections of $Pz$ and $Q$, so that $Pz=\bigoplus_{i=1}^m\mathbb Cp_i$ and $Q=\bigoplus_{j=1}^n\mathbb Cq_j$.
Since $P\subset_{\delta}Q$, Lemma \ref{Ch} provides a $*$-homomorphism $\theta:P\rightarrow\mathbb M_2(\mathbb C)\otimes Q$ such that $\|\theta(x)-e_{1,1}\otimes x\|_{2,\widetilde\tau}\leq 30\delta^{\frac{1}{8}}$, for every $x\in (P)_1$. In particular, using \eqref{zzz} we get that \begin{equation}\label{z}\|\theta(z)-e_{1,1}\otimes 1\|_{2,\widetilde\tau}\leq \|\theta(z)-e_{1,1}\otimes z\|_{2,\widetilde\tau}+\|e_{1,1}\otimes (1-z)\|_{2,\widetilde\tau}\leq 30\delta^{\frac{1}{8}}+5\delta\leq 31\delta^{\frac{1}{8}}.\end{equation}

Write $\theta(z)=\sum_{j=1}^n\alpha_j\otimes q_j$, where $\alpha_j\in\mathbb M_2(\mathbb C)$ is a projection, for every $1\leq j\leq n$. Let $S$ be the set of all $j\in\{1,\cdots,n\}$ such that $\alpha_j$ has rank one. Define $w=\sum_{j\in S}q_j\in Q$.
If $j\notin S$, then $\alpha_j$ is equal to $0$ or $1$ and thus  $\|\alpha_j\otimes q_j-e_{1,1}\otimes q_j\|_{2,\widetilde\tau}^2=\frac{\tau(q_j)}{2}$.  This implies that $$\|1-w\|_{2,\tau}^2=\sum_{j\notin S}\tau(q_j)=2\sum_{j\notin S}\|\alpha_j\otimes q_j-e_{1,1}\otimes q_j\|_{2,\widetilde\tau}^2\leq 2\|\theta(z)-e_{1,1}\otimes 1\|_{2,\widetilde\tau}^2.$$ In combination with \eqref{z}, we derive that \begin{equation}\label{w}\|1-w\|_{2,\tau}\leq 31\sqrt{2}\delta^{\frac{1}{8}}.\end{equation}

Since $\theta(z)(1\otimes w)=\sum_{j\in S}\alpha_j\otimes q_j$ and $\alpha_j$ has rank one, for every $j\in S$, we get that 
there is a partition $S=S_1\sqcup\cdots\sqcup S_m$ such that $\theta(p_i)(1\otimes w)=\sum_{j\in S_i}\alpha_j\otimes q_j$, for every $1\leq i\leq m$. 
Define a $*$-homomorphism $\rho:P\rightarrow Q$ be letting $\rho(1-z)=0$ and $\rho(p_i)=\sum_{j\in S_i}q_j$, for every $1\leq i\leq m$. 

\begin{claim}\label{rho}
$\|\rho(x)-x\|_{2,\tau}\leq 150\delta^{\frac{1}{8}}$, for every $x\in (P)_1$.
\end{claim}

{\it Proof of Claim \ref{rho}.}
Let $x\in (P)_1$ and write $x=c_0+\sum_{i=1}^mc_ip_i$, where $c_0\in (P(1-z))_1$ and $c_i\in \mathbb C$ satisfies $|c_i|\leq 1$, for $0\leq i\leq m$. Let $y=\sum_{i=1}^mc_ip_i$. Then $e_{1,1}\otimes \rho(x)=\sum_{i=1}^m\sum_{j\in S_i}c_i(e_{1,1}\otimes q_j)$ and $\theta(y)(1\otimes w)=\sum_{i=1}^mc_i\theta(p_i)(1\otimes w)=\sum_{i=1}^m\sum_{j\in S_i}c_i(\alpha_j\otimes q_j)$. Since the projections $\{q_j\}_{j\in S}$ are pairwise orthogonal, the sets $\{S_i\}_{i=1}^m$ partition $S$ and $|c_i|\leq 1$, for every $1\leq i\leq m$, we get that $$\|e_{1,1}\otimes \rho(x)-\theta(y)(1\otimes w)\|_{2,\widetilde\tau}^2\leq\sum_{j\in S}\|e_{1,1}\otimes q_j-\alpha_j\otimes q_j\|_{2,\widetilde\tau}^2\leq \|e_{1,1}\otimes 1-\theta(z)\|_{2,\widetilde\tau}^2.$$ In combination with \eqref{z}, we derive that \begin{equation}\label{rh} \|e_{1,1}\otimes \rho(x)-\theta(y)(1\otimes w)\|_{2,\widetilde\tau}\leq 31\delta^{\frac{1}{8}}.\end{equation}

Since $\|\theta(y)-e_{1,1}\otimes y\|_{2,\widetilde\tau}\leq 30\delta^{\frac{1}{8}}$, using \eqref{w} and \eqref{rh}, we further get that 
\begin{align*}\|\rho(x)-y\|_{2,\tau}&=\sqrt{2}\|e_{1,1}\otimes\rho(x)-e_{1,1}\otimes y\|_{2,\widetilde\tau}\\&\leq \sqrt{2}\big(\|e_{1,1}\otimes\rho(x)-\theta(y)(1\otimes w)\|_{2,\widetilde\tau}+\|\theta(y)(1\otimes (1-w))\|_{2,\widetilde\tau}+\|\theta(y)-e_{1,1}\otimes y\|_{2,\widetilde\tau}\big)\\&\leq 140\delta^{\frac{1}{8}}.\end{align*}

As $x-y\in (P(1-z))_1$, \eqref{zz} implies that $\|x-y\|_{2,\tau}\leq \|1-z\|_{2,\tau}\leq 5\delta$. The last displayed inequality gives that $\|\rho(x)-x\|_{2,\tau}\leq \|\rho(x)-y\|_{2,\tau}+\|1-z\|_{2,\tau}\leq 140\delta^{\frac{1}{8}}+5\delta\leq 150\delta^{\frac{1}{8}}$, proving the claim. \hfill$\square$

Finally, Claim \ref{rho} gives that ${\bf d}_\tau(P,\rho(P))\leq 150\delta^{\frac{1}{8}}$. Let $R=\rho(P)\bigoplus\mathbb C(1-w)$. Since $\rho(1)=w$, $R$ is  von Neumann subalgebra of $Q$. By \eqref{w} we get that ${\bf d}_\tau(\rho(P),R)\leq \|1-w\|_{2,\tau}\leq 31\sqrt{2}\delta^{\frac{1}{8}}\leq 50\delta^{\frac{1}{8}}$. Thus, we conclude that ${\bf d}_\tau(P,R)\leq {\bf d}_\tau(P,\rho(P))+{\bf d}_\tau(\rho(P),R)\leq 200\delta^{\frac{1}{8}}.$
\hfill$\blacksquare$

The following consequence of Corollary \ref{Chr} is a key ingredient of the proof of Proposition \ref{reduction}.

\begin{lemma}\label{perturb} For any $\varepsilon>0$, there is $\delta=\delta_3(\varepsilon)>0$ such that the following holds.
For $i\in\{1,2,3\}$, let $M_i$ be a finite dimensional factor and denote by $1_{M_i}$ its unit.
 Let $M=M_1\otimes M_2\otimes M_3$ 
 and $P\subset M$ be a von Neumann subalgebra. Assume that $M_1\otimes 1_{M_2}\otimes 1_{M_3}\subset_{\delta}P$ and $P\subset_{\delta}M_1\otimes M_2\otimes 1_{M_3}$.
Then there exists a von Neumann subalgebra  $S\subset M_2$ such that 
$\emph{\bf{d}}(P,M_1\otimes S\otimes 1_{M_3})\leq \varepsilon$.

\end{lemma}

{\it Proof of Lemma \ref{perturb}.} Let $\delta_1:(0,+\infty)\rightarrow (0,+\infty)$ be the function provided by Corollary \ref{Chr}. 
Let $\varepsilon>0$. Let $\varepsilon'>0$ such that $\varepsilon'\leq\frac{\varepsilon}{2}$ and $\delta(\varepsilon')+\varepsilon'\leq\frac{1}{4}\delta_1(\frac{\varepsilon}{8})$. 
We will show that $\delta:=\delta(\varepsilon')$ works.

Since $M_1,M_2$ are factors, so is $M_1\otimes M_2$. Since $P\subset_{\delta(\varepsilon')}M_1\otimes M_2\otimes 1_{M_3}$, Corollary \ref{Chr} implies the existence of a von Neumann subalgebra $Q\subset M_1\otimes M_2$ such that 
\begin{equation}
\label{unu}\text{\bf{d}}(P,Q\otimes 1_{M_3})\leq \varepsilon'\leq\frac{\varepsilon}{2}.
\end{equation}

Since $M_1\otimes 1_{M_2}\otimes 1_{M_3}\subset_{\delta_1(\varepsilon')}P$,  \eqref{unu} implies that $M_1\otimes 1_{M_2}\subset_{\delta_1(\varepsilon')+\varepsilon'}Q$.  Let $N=Q'\cap (M_1\otimes M_2)$.  Since $(M_1\otimes 1_{M_2})'\cap (M_1\otimes M_2)=1_{M_1}\otimes M_2$ and $4(\delta_1(\varepsilon')+\varepsilon')\leq\delta_1(\frac{\varepsilon}{8})$, Lemma \ref{commutant} gives that \begin{equation}\label{doi}N\subset_{\delta_1(\frac{\varepsilon}{8})}1_{M_1}\otimes M_2.\end{equation}

Since $M_2$ is a factor, applying Corollary \ref{Chr} gives a von  Neumann subalgebra $R\subset M_2$ such that
\begin{equation}\label{trei}
\text{\bf{d}}(N,1_{M_1}\otimes R)\leq \frac{\varepsilon}{8}.
\end{equation}
As $M_1\otimes M_2$ is a finite dimensional factor, the bicommutant theorem gives that $N'\cap (M_1\otimes M_2)=Q$. By applying Lemma \ref{commutant} again we deduce that if $S=R'\cap M_2$, then 
\begin{equation}
\label{patru}
\text{\bf{d}}(Q,M_1\otimes S)\leq \frac{\varepsilon}{2}.
\end{equation}
Finally, by combining \eqref{unu} and \eqref{patru} we derive that $\text{\bf{d}}(P,M_1\otimes S\otimes 1_{M_3})\leq \varepsilon$.
 \hfill$\blacksquare$

\section{Pairs of unitary matrices with spectral gap}

 The goal of this section is to prove the following two results giving pairs
 of unitary matrices with spectral gap properties. These provide the first step towards proving Theorem \ref{main}. For $n\in\mathbb N$, we denote by $\tau$ the normalized trace on $\mathbb M_n(\mathbb C)$, and by $\|\cdot\|_2$ and $\|\cdot\|_1$ the associated norms.

\begin{proposition}\label{spF3} There is a constant $\eta>0$ such that the following holds. Let $A=\mathbb M_k(\mathbb C)$ and $B=\mathbb M_n(\mathbb C)$, for $k,n\in\mathbb N$. 
Then there are $Z_1,Z_2\in\mathcal U(A\otimes 1)$
such that $$\text{$\|x-\emph{E}_{1\otimes B}(x)\|_2\leq \eta(||[Z_1,x]||_2+||[Z_2,x]||_2)$, for every $x\in A\otimes B$.}$$ 
\end{proposition}

Proposition \ref{spF3} sufficces to prove that $\mathbb F_3\times\mathbb F_3$, and thus $\mathbb F_m\times\mathbb F_n$, for all $m,n\geq 3$, is not HS-stable. However, to prove the failure of HS-stability for $\mathbb F_2\times\mathbb F_2$, we will need the following result. 

\begin{proposition}\label{spF2} 
There is a constant $\eta>0$ such that the following holds. Let $A=\mathbb M_k(\mathbb C)$, $B=\mathbb M_n(\mathbb C)$ and $w\in \mathcal U(A\otimes B)$, for $k,n\in\mathbb N$. 
Then there are $Z_1,Z_2\in\mathcal U(\mathbb M_3(\mathbb C)\otimes A\otimes B)$
such that 
 \begin{enumerate}
 \item $\|x-\emph{E}_{1\otimes 1\otimes B}(x)\|_2\leq \eta(||[Z_1,x]||_2+||[Z_2,x]||_2)$, for every $x\in \mathbb M_3(\mathbb C)\otimes A\otimes B$, 
 \item $Z_1\in \mathbb M_3(\mathbb C)\otimes A\otimes 1$, and 
 \item  $Z_2=\begin{pmatrix}0&0&1\\1&0&0\\0&w&0\end{pmatrix}$.
 \end{enumerate}
\end{proposition}

\subsection{Pairs of unitary matrices with spectral gap} 
The proofs of Propositions \ref{spF3} and \ref{spF2} rely  on the following result.

\begin{lemma}\label{gap}
There exist a constant $\kappa>0$, a sequence $(k_n)$ of natural numbers with $k_n\rightarrow\infty$, and a pair of unitaries $(u_n,v_n)\in \emph{U}(k_n)^2$, for every $n\in\mathbb N$, such that $$\text{$\|x-\tau(x)1\|_2\leq \kappa(\|[u_n,x]\|_2+\|[v_n,x]\|_2)$, for every $x\in\mathbb M_{k_n}(\mathbb C)$.}$$

Moreover, we can take 
$k_n=n$, for every $n\in\mathbb N$. 
\end{lemma}

This result is likely known to experts, but, for completeness, we indicate how it follows from the literature.
We give two proofs of the main assertion based on property (T) and quantum expanders, respectively. The second proof will allow us to also derive the moreover assertion. 

\subsection*{First proof of the main assertion of Lemma \ref{gap}} The first proof combines an argument from the proof of \cite[Proposition 3.9(4)]{BCI15}, which we recall below, 
with the fact that $\Gamma:=\text{SL}_3(\mathbb Z)$ is $2$-generated. Indeed, by \cite{Tr62}, the following two matrices generate $\Gamma$: $$\text{$a=\begin{pmatrix}0&1&0\\0&0&1\\1&0&0\end{pmatrix}$ and $b=\begin{pmatrix}1&0&0\\1&1&0\\0&0&1 \end{pmatrix}$.}$$ 
Since $\Gamma$ has Kazhdan's property (T) (see, e.g., \cite[Theorem 12.1.14]{BO08}), we can find $\kappa>0$ such that if $\rho:\Gamma\rightarrow \text{U}(\mathcal H)$ is any unitary representation and $P:\mathcal H\rightarrow\mathcal H$ is the orthogonal projection onto the subspace of $\rho(\Gamma)$-invariant vectors, then \begin{equation}\label{TT}\text{$\|\xi-P(\xi)\|\leq\kappa (\|\rho(a)\xi-\xi\|+\|\rho(b)\xi-\xi\|)$, for every $\xi\in\mathcal H$.}\end{equation}

Since $\Gamma$ is residually finite and has property (T), it has a sequence of finite dimensional irreducible representations $\pi_n:\Gamma\rightarrow \text{U}(k_n)$, $n\in\mathbb N$, with $k_n\rightarrow\infty$ (see the proof of \cite[Proposition 3.9(4)]{BCI15}).
Alternatively, if $p$ is a prime, then any nontrivial representation of $\text{SL}_3(\mathbb Z/p\mathbb Z)$  has dimension at least $\frac{p-1}{2}$ (see \cite[Exercise 3.0.9]{Ta15}). Thus, we can take $\pi_n$ to be any irreducible representation of $\Gamma$  factoring through $\text{SL}_3(\mathbb Z/p_n\mathbb Z)$, for any $n\in\mathbb N$, where $(p_n)$ is a sequence of primes with $p_n\rightarrow\infty$.

Since $\pi_n$ is irreducible, the only matrices which are invariant under the unitary representation $\rho_n:\Gamma\rightarrow \text{U}(\mathbb M_{k_n}(\mathbb C))$ 
given by $\rho_n(g)x=\pi_n(g)x\pi_n(g)^*$ are the scalar multiples of the identity. 
Thus, applying inequality \eqref{TT} to $\rho_n$ gives that $\|x-\tau(x)1\|_2\leq\kappa (\|[\pi_n(a),x]\|_2+\|[\pi_n(b),x]\|_2)$, for every $x\in\mathbb M_{k_n}(\mathbb C)$.  Hence, $u_n=\pi_n(a)$ and $v_n=\pi_n(b)$ satisfy the conclusion of Lemma \ref{gap}.
\hfill$\blacksquare$

We next give a second proof of the main assertion of Lemma \ref{gap} showing that one can take $k_n=n$.
This relies on the notion on quantum expanders introduced in \cite{BS07,Ha07} (see also \cite{Pi12}).
For a related application of quantum expanders, see the recent article \cite{MS20}.
For $k\geq 2$ and a $k$-tuple of unitaries $u=(u_1,u_2,...,u_k)\in \text{U}(n)^k$,
let $T_u:\mathbb M_n(\mathbb C)\rightarrow\mathbb M_n(\mathbb C)$ be the operator given by $$\text{$T_u(x)=\sum_{i=1}^ku_ixu_i^*$, for every $x\in\mathbb M_n(\mathbb C)$}.$$
Endow $\mathbb M_n(\mathbb C)$ with the normalized Hilbert-Schmidt norm, note that the space $\mathbb M_n(\mathbb C)\ominus\mathbb C1$ of matrices of trace zero is $T_u$-invariant, and denote by $T_u^0$ the restriction of $T_u$ to $\mathbb M_n(\mathbb C)\ominus\mathbb C1$. 

\begin{remark}\label{k=2} We clearly have that $\|T_u^0\|\leq k$. Moreover, equality holds if $k=2$. To see this, let $u=(u_1,u_2)$. Then $T_u(u_2^*u_1)=2u_1u_2^*$ and $\|u_2^*u_1-\alpha 1\|_2=\|u_1u_2^*-\alpha 1\|_2=\sqrt{1-|\alpha|^2}$, where $\alpha=\tau(u_2^*u_1)=\tau(u_1u_2^*)$. If $|\alpha|<1$, then since $T_u(u_2^*u_1-\alpha 1)=2(u_1u_2^*-\alpha 1)$ we get that $\|T_u^0\|=2$. If $|\alpha|=1$, then $u_1=\alpha u_2$ and so $T_u(x)=2u_1xu_1^*$, for every $x\in\mathbb M_n(\mathbb C)$, which gives that $\|T_u^0\|=2$.
\end{remark}

A sequence of $k$-tuples $u^n=(u_1^n,...,u_k^n)\in \text{U}(n)^k$ is called a {\it quantum expander} if $\sup_n\|T_{u^n}^0\|<k$. By Remark \ref{k=2}, this forces that $k\geq 3$.
The following result due to Hastings \cite{Ha07} (formulated here following \cite[Lemma 1.8]{Pi12}) shows that random unitaries provide quantum expanders for $k\geq 3$.
\begin{lemma}\label{hastings} For $n\in\mathbb N$, let $\mu_n$ be the Haar measure of $\emph{U}(n)$. Then for every $\varepsilon>0$ we have that $$\lim\limits_{n\rightarrow\infty}\mu_n^k(\{u\in \emph{U}(n)^k\mid  \|T_u^0\|\leq 2\sqrt{k-1}+\varepsilon k \})=1.$$

\end{lemma}

\subsection*{Proof of Lemma \ref{gap}} We claim there is $N\in\mathbb N$ such that the main assertion of Lemma \ref{gap} holds for any constant $\kappa$ greater than $3+2\sqrt{2}$ and $k_n=n$, for every $n>N$. Assuming this claim, note that if $n\in\mathbb N$ is fixed, then we can find two unitaries $u_n,v_n\in\text{U}(n)$ such that $\{u_n,v_n\}'\cap \mathbb M_n(\mathbb C)=\mathbb C1$. Using the compactness of the unit ball of $\mathbb M_n(\mathbb C)$ with respect to the $\|\cdot\|_2$-norm, we can find a constant $\kappa_n>0$ such that $\|x-\tau(x)1\|_2\leq\kappa_n(\|[u_n,x]\|_2+\|[v_n,x]\|_2)$, for every $x\in\mathbb M_n(\mathbb C)$. 
It is now clear that the moreover assertion of Lemma \ref{gap} holds after replacing $\kappa$ with $\max\{\kappa,\kappa_1,...,\kappa_N\}$.

To prove our claim, fix a constant $\kappa>3+2\sqrt{2}$. Note that $\frac{1}{\kappa}<\frac{1}{3+2\sqrt{2}}=3-2\sqrt{2}$ and let $\varepsilon:=\frac{3-2\sqrt{2}-\frac{1}{\kappa}}{2}>0$.
By applying Lemma \ref{hastings} in the case $k=3$ we deduce that \begin{equation}\label{k=3}\lim\limits_{n\rightarrow\infty}\mu_n^3(\{u\in \text{U}(n)^3\mid \|T_u^0\|\leq 2\sqrt{2}+2\varepsilon\})=1.\end{equation} Let $S_n$ be set of pairs of unitaries $(u_1,u_2)\in \text{U}(n)^2$ such that $\|T_{(u_1,u_2,I)}^0\|\leq 2\sqrt{2}+2\varepsilon$. Since $||T_{(u_1,u_2,u_3)}^0||=||T_{(u_3^*u_1,u_3^*u_1,I)}^0||$, for every $(u_1,u_2,u_3)\in \text{U}(n)^3$, \eqref{k=3} implies that $\lim\limits_{n\rightarrow\infty}\mu_n^2(S_n)=1$. 

Now, let $(u_1,u_2)\in S_n$. Then for every $x\in\mathbb M_n(\mathbb C)\ominus\mathbb C1$ we have \begin{align*}(2\sqrt{2}+2\varepsilon)\|x\|_2\geq \|u_1xu_1^*+u_2xu_2^*+x\|_2\geq 3\|x\|_2-\|[u_1,x]\|_2-\|[u_2,x]\|_2,
\end{align*} and hence $\|x\|_2\leq \kappa(\|[u_1,x]\|_2+\|[u_2,x]\|_2)$. If $x\in\mathbb M_n(\mathbb C)$, then applying this last inequality to $x-\tau(x)1\in\mathbb M_n(\mathbb C)\ominus\mathbb C1$ gives that $\|x-\tau(x)1\|_2\leq \kappa(\|[u_1,x]\|_2+\|[u_2,x]\|_2)$.

Since $\lim\limits_{n\rightarrow\infty}\mu_n^2(S_n)=1$, we have that $S_n\not=\emptyset$, for $n$ large enough. Then any pair $(u_n,v_n)\in S_n$, for $n$ large enough, will satisfy the conclusion of Lemma \ref{gap}.
\hfill$\blacksquare$

We end this subsection by proving the following property for pairs of unitaries with spectral gap:
\begin{lemma}\label{corner}
Let  $\kappa>0$ and $(u_1,u_2)\in \emph{U}(n)^2$ such that $\|x-\tau(x)1\|_2\leq\kappa(\|[u_1,x]\|_2+\|[u_2,x]\|_2)$, for all $x\in\mathbb M_n(\mathbb C)$. Then  $\|x\|_2\leq 10^5\kappa^6(\|u_1xv-x\|_2+\|u_2xv-x\|_2)$, for all $v\in \emph{U}(n)$ and $x\in\mathbb M_n(\mathbb C)$.
\end{lemma}

{\it Proof.} 
We claim that $\delta:=\|u_1-u_2\|_2\geq \frac{1}{8\kappa^2}$. Indeed, $\|u_1-\tau(u_1)1\|_2\leq \kappa\|[u_1,u_2]\|_2\leq 2\kappa\delta$ and similarly $\|u_2-\tau(u_2)1\|_2\leq 2\kappa\delta$. If $x\in \text{U}(n)$ and $\tau(x)=0$, then $\|[u_1,x]\|_2\leq 2\|u_1-\tau(u_1)1\|_2\leq 4\kappa\delta$ and similarly $\|[u_2,x]\|_4\leq 2\kappa\delta$.
Thus,
$1=\|x\|_2\leq \kappa(\|[u_1,x]\|_2+\|[u_2,x]\|_2)\leq 8\kappa^2\delta$, proving our claim.
Note also that since $1=\|x\|_2\leq \kappa(\|[u_1,x]\|_2+\|[u_2,x]\|_2)\leq 4\kappa$, we have that $\kappa\geq\frac{1}{4}$.

To prove the conclusion of the lemma, let $x\not=0$ and put $\varepsilon=\|x\|_2^{-1}(\|u_1xv-x\|_2+\|u_2xv-x\|_2)$. Let $w\in \{u_1,u_2\}$. Then $wxx^*w^*=(wxv)(wxv)^*$ and the Cauchy-Schwarz inequality implies that 
$$
\|wxx^*w^*-xx^*\|_1\leq \|(wxv-x)(wxv)^*\|_1+\|x(wxv-x)^*\|_1\leq 2\varepsilon \|x\|_2^2.
$$
Let $y=(xx^*)^{\frac{1}{2}}$ and $\alpha=\tau(y)\geq 0$. Then the Powers-St{\o}rmer inequality \eqref{PS} implies that 
$$\text{$\|wyw^*-y\|_2\leq \|wy^2w^*-y^2\|_1^{\frac{1}{2}}=\|wxx^*w^*-xx^*\|_1^{\frac{1}{2}}
\leq\sqrt{2\varepsilon}\|x\|_2$, for every $w\in\{u_1,u_2\}.$}$$

Thus, we get that $\|y-\alpha 1\|_2\leq 2\kappa\sqrt{2\varepsilon}\|x\|_2$. Since  $\|y\|_2=\|x\|_2$, it follows that \begin{equation}\label{alpha}\alpha\geq (1-2\kappa\sqrt{2\varepsilon})\|x\|_2.\end{equation}
 
Next, by the polar decomposition we can find a unitary $z\in \text{U}(n)$ such that $x=yz$. Then we have $\|x-\alpha z\|_2=\|y-\alpha 1\|_2\leq 2\kappa\sqrt{2\varepsilon}\|x\|_2$. Since $\|(u_1-u_2)xv\|_2\leq \|u_1xv-x\|_2+\|u_2xv-x\|_2=\varepsilon\|x\|_2$, we further get that 
$\alpha\delta=\|(u_1-u_2)(\alpha z)v\|_2\leq \|(u_1-u_2)xv\|_2+\|x-\alpha z\|_2\leq (\varepsilon+2\kappa\sqrt{2\varepsilon})\|x\|_2.$

Since $x\not=0$, by combining the fact that $\delta\geq \frac{1}{8\kappa^2}$ with the last inequality we conclude that $$\frac{1-2\kappa\sqrt{2\varepsilon}}{8\kappa^2}\leq\varepsilon+2\kappa\sqrt{2\varepsilon}.$$

Since $\kappa\geq\frac{1}{4}$ and $\varepsilon\leq 4$, it follows that $\varepsilon\geq\frac{1}{10^5\kappa^6}$, which finishes the proof.
\hfill$\blacksquare$

\subsection{Proof of Proposition \ref{spF3}}  Let $\kappa>0$ be as given by Lemma \ref{gap}. The moreover assertion of Lemma \ref{gap} provides $u,v\in\mathcal U(A)$
 such that $\|x-\tau(x)\|_2\leq\kappa(\|[u,x]\|_2+\|[v,x]\|_2)$, for every $x\in A$. 
 
It is a standard fact that $Z_1=u\otimes 1, Z_2=v\otimes 1\in\mathcal U(A\otimes 1)$ satisfy the conclusion for $\eta=\sqrt{2}\kappa$. For completeness, let us recall the argument. Let $\{\xi_i\}_{i\in I}$ be an orthonormal basis of $B$ with respect to the scalar product given by its trace. Let $x\in A\otimes B$ and write $x=\sum_{i\in I}x_i\otimes \xi_i$, with $x_i\in A$. Then $\text{E}_{1\otimes B}(x)=\sum_{i\in I}\tau(x_i)1\otimes \xi_i, [Z_1,x]=\sum_{i\in I}[u,x_i]\otimes\xi_i, [Z_2,x]=\sum_{i\in I}[v,x_i]\otimes\xi_i$ and therefore 
\begin{align*}
\|x-\text{E}_{1\otimes B}(x)\|_2^2=\sum_{i\in I}\|x_i-\tau(x_i)1\|_2^2&\leq\sum_{i\in I}\kappa^2\big(\|[u,x_i]\|_2+\|[v,x_i]\|_2\big)^2
\\&\leq \eta^2\sum_{i\in I}\big(\|[u,x_i]\|_2^2+\|[v,x_i]\|_2^2\big)\\&=\eta^2\big(\|[Z_1,x]\|_2^2+\|[Z_2,x]\|_2^2\big)\\&\leq \Big(\eta\big(\|[Z_1,x]\|_2+\|[Z_2,x]\|_2\big)\Big)^2.
\end{align*}

This finishes the proof. \hfill$\blacksquare$

\subsection{Proof of Proposition \ref{spF2}} Let $\kappa>0$ be as given by Lemma \ref{gap}. The moreover assertion of Lemma \ref{gap} gives $u,v\in\mathcal U(A)$ such that $\|x-\tau(x)1\|_2\leq \kappa(\|[u,x]\|_2+\|[v,x]\|_2)$, for every $x\in A$.  The proof of Proposition \ref{spF3} shows that
\begin{equation}\label{111}
\text{$\|x-\text{E}_{1\otimes B}(x)\|_2\leq \sqrt{2}\kappa(\|[u\otimes 1,x]\|_2+\|[v\otimes 1,x]\|_2)$, for every $x\in A\otimes B$.}
\end{equation}
By Lemma \ref{corner}, $\|x\|_2\leq 10^5\kappa^6(\|uxt-x\|_2+\|vxt-x\|_2)$, for every $t\in\mathcal U(A)$ and $x\in A$. Then the argument from the proof of Proposition \ref{spF3} implies that for every $t\in\mathcal U(A)$ and $x\in A\otimes B$ we have \begin{equation}\label{222}\text{$\|x\|_2\leq \sqrt{2}\cdot10^5\kappa^6(\|(u\otimes 1)x(t\otimes 1)-x\|_2+\|(v\otimes 1)x(t\otimes 1)-x\|_2)$.} \end{equation}

Define $Z_1,Z_2\in\mathcal U(\mathbb M_3(\mathbb C)\otimes A\otimes B)$ by letting $$\text{$Z_1=\begin{pmatrix} u\otimes 1&0&0\\0& u\otimes 1 &0\\0&0&v\otimes 1\end{pmatrix}$ \;\; and \;\; $Z_2=\begin{pmatrix}0&0&1\\1&0&0\\0&w&0\end{pmatrix}$.} $$
Then $Z_1,Z_2$  satisfy conditions (2) and (3) from the conclusion and
$$Z_2Z_1Z_2^*=\begin{pmatrix}v\otimes 1&0&0\\0&u\otimes 1&0\\0&0& w(u\otimes 1)w^* \end{pmatrix}.$$

We will show that condition (1) is satisfied for $\eta=10^7(1+\kappa^6)$. 
To this end, fix $x\in\mathbb M_3(\mathbb C)\otimes A\otimes B$ with $x=x^*$. Write $x=[x_{ij}]$, where $x_{i,j}\in A\otimes B$ are such that $x_{i,j}^*=x_{j,i}$, for every $1\leq i,j\leq 3$. Our goal is to show that 
\begin{equation}\label{goal}
\|x-\text{E}_{1\otimes 1\otimes B}(x)\|_2\leq\frac{\eta}{2}(\|[Z_1,x]\|_2+\|[Z_2,x]\|_2).
\end{equation}

Towards this goal, we denote $\varepsilon=\|[Z_1,x]\|_2+\|[Z_2,x]\|_2$ and record the following elementary fact:

\begin{fact}\label{fact}
Let $d_1,d_2,d_3\in \mathcal U(A\otimes B)$ and put $d=[d_i\delta_{i,j}]\in\mathcal U(\mathbb M_3(\mathbb C)\otimes A\otimes B)$. Then we have $\|[d,x]\|_2^2
=\sum_{i,j=1}^3\|d_ix_{i,j}d_j^*-x_{i,j}\|_2^2$, hence $\|d_ix_{i,j}d_j^*-x_{i,j}\|_2\leq \|[d,x]\|_2$, for every $1\leq i,j\leq 3$.
\end{fact}

Fact \ref{fact} implies that $\|[(u\otimes 1),x_{1,1}]\|_2\leq \|[Z_1,x]\|_2\leq\varepsilon$ and $\|[(v\otimes 1),x_{1,1}]\|_2\leq \|[Z_2Z_1Z_2^*,x]\|_2\leq 3\varepsilon$.
Together with \eqref{111} this gives that
\begin{equation}\label{un}
\|x_{1,1}-y\|_2\leq (\sqrt{2}\kappa)\cdot (4\varepsilon)\leq\frac{\eta}{6}\varepsilon,\end{equation}
where $y:=\text{E}_{1\otimes B}(x_{1,1})\in 1\otimes B$.

By using Fact \ref{fact}, we also derive that 
$\|(u\otimes 1)x_{1,2}(u\otimes 1)^*-x_{1,2}\|_2\leq \|[Z_1,x]\|_2\leq \varepsilon$ and that $\|(v\otimes 1)x_{1,2}(u\otimes 1)^*-x_{1,2}\|_2\leq \|[Z_2Z_1Z_2^*,x]\|_2\leq 3\varepsilon$.
Applying \eqref{222} to $x=x_{1,2}$ and $t=u^*$, we get
\begin{equation}\label{trois}\|x_{1,2}\|_2\leq (\sqrt{2}\cdot 10^5\kappa^6)\cdot (4\varepsilon)\leq\frac{\eta}{6}\varepsilon.
\end{equation}
Next, note that 
$$[Z_2,x]=\begin{pmatrix}x_{3,1}-x_{1,2}&x_{3,2}-x_{1,3}w&x_{3,3}-x_{1,1}\\ x_{1,1}-x_{2,2}&x_{1,2}-x_{2,3}w&x_{1,3}-x_{2,1}\\wx_{2,1}-x_{3,2}&wx_{2,2}-x_{3,3}w&wx_{2,3}-x_{3,1} \end{pmatrix}$$
Since $\|[Z_2,x]\|_2\leq\varepsilon$, we get that $\|x_{2,2}-x_{1,1}\|_2\leq \varepsilon$ and $\|x_{3,3}-x_{1,1}\|_2\leq\varepsilon$. Together with \eqref{un}, this gives that 
\begin{equation}\label{dois} \text{$\|x_{2,2}-y\|_2\leq (\sqrt{2}\cdot 8\kappa+1)\varepsilon\leq\frac{\eta}{6}\varepsilon$ and $\|x_{3,3}-y\|_2\leq (\sqrt{2}\cdot 8\kappa+1)\varepsilon\leq\frac{\eta}{6}\varepsilon$.}
\end{equation}
Since $\|[Z_2,x]\|_2\leq\varepsilon$, we also get that  $\|x_{2,3}w-x_{1,2}\|_2\leq\varepsilon$ and $\|x_{3,1}-x_{1,2}\|_2\leq\varepsilon$. Together with \eqref{trois}, this gives that 
\begin{equation}\label{quatre}\text{$\|x_{2,3}\|_2\leq (8\cdot\sqrt{2}\cdot 10^4\kappa^6+1)\varepsilon\leq\frac{\eta}{6}\varepsilon$ \;\; and \;\;$\|x_{3,1}\|_2\leq (8\cdot\sqrt{2}\cdot 10^4\kappa^6+1)\varepsilon\leq\frac{\eta}{6}\varepsilon$.}
\end{equation}
Since $x=x^*$, by using \eqref{un}, \eqref{trois}, \eqref{dois} and \eqref{quatre}, we get that
\begin{align*}
\|x-1\otimes y\|_2^2&=\|x_{1,1}-y\|_2^2+\|x_{2,2}-y\|_2^2+\|x_{3,3}-y\|_2^2+2\|x_{1,2}\|_2^2+2\|x_{1,3}\|_2^2+2\|x_{2,3}\|_2^2\leq 9\big(\frac{\eta}{6}\varepsilon\big)^2.
\end{align*}
Since $1\otimes y\in 1\otimes 1\otimes B$, we get that $\|x-\text{E}_{1\otimes 1\otimes B}(x)\|_2\leq \|x-1\otimes y\|_2\leq\frac{\eta}{2}\varepsilon$, hence \eqref{goal} holds. 

Finally, given $x\in \mathbb M_3(\mathbb C)\otimes A\otimes B$, write $x=x_1+ix_2$, where $x_1=x_1^*$ and $x_2=x_2^*$. 
Then $\|[u,x]\|_2^2=\|[u,x_1]\|_2^2+\|[u,x_2]\|_2^2$, for every unitary $u$ and by using \eqref{goal} for $x_1$ and $x_2$ we get that
\begin{align*}
\|x-\text{E}_{1\otimes 1\otimes B}(x)\|_2&\leq \|x_1-\text{E}_{1\otimes 1\otimes B}(x_1)\|_2+\|x_2-\text{E}_{1\otimes 1\otimes B}(x_2)\|_2 \\&\leq\frac{\eta}{2}(\|[Z_1,x_1]\|_2+\|[Z_2,x_1]\|_2+\|[Z_1,x_2]\|_2+\|[Z_2,x_2]\|_2)\\&\leq \eta(\|[Z_1,x]\|_2+\|[Z_2,x]\|_2).
\end{align*} This finishes the proof.
\hfill$\blacksquare$

\section{Proof of Proposition \ref{reduction}}
This section is devoted to the proof of Proposition \ref{reduction}. 
We first prove Proposition \ref{reduction} under the stronger assumption that $\mathbb F_3\times\mathbb F_3$ is \text{HS}-stable, since the proof is more transparent in this case and relies on the simpler Proposition \ref{spF3} instead of Proposition \ref{spF2}.

\subsection*{Proof of Proposition \ref{reduction} assuming that $\mathbb F_3\times\mathbb F_3$ is \text{HS} stable}
Let $\varepsilon>0$. 
Let $\eta>0$ be the constant provided by Proposition \ref{spF3}. 
Let $\delta_3:(0,+\infty)\rightarrow (0,+\infty)$ be the function provided by Lemma \ref{perturb}. Let $\varepsilon_0>0$ such that $\varepsilon_0<\frac{\varepsilon}{24}$ and $16\eta\varepsilon_0<\delta_3(\frac{\varepsilon}{16})$.

Since $\mathbb F_3\times\mathbb F_3$ is \text{HS}-stable we can find $\delta>0$ such that for any finite dimensional factor $M$ and $Z_\alpha,T_\beta\in\mathcal U(M)$ such that $\|[Z_\alpha,T_\beta]\|_2^2\leq\delta$, for every $\alpha,\beta\in\{1,2,3\}$, we can find $\widetilde Z_\alpha,\widetilde T_\beta\in\mathcal U(M)$ such that $\|\widetilde Z_\alpha-Z_\alpha\|_2\leq\varepsilon_0, \|\widetilde T_\beta-T_\beta\|_2\leq\varepsilon_0$ and $[\widetilde Z_\alpha,\widetilde T_\beta]=0$, for every $\alpha,\beta\in\{1,2,3\}$.

Let $k,m,n\in\mathbb N$ and $U_1,...,U_k,V_1,...,V_m\in U(n)$ such that \begin{equation}\label{almost}\frac{1}{km}\sum_{i=1}^k\sum_{j=1}^m\|[U_i,V_j]\|_2^2\leq\delta.\end{equation} 
Denote $A=\mathbb M_k(\mathbb C), B=\mathbb M_n(\mathbb C)$, $C=\mathbb M_m(\mathbb C)$ and $M=A\otimes B\otimes C$. By applying Lemma \ref{spF3} twice, we can find $Z_1,Z_2\in\mathcal U(A\otimes 1\otimes 1)$ and $T_1,T_2\in\mathcal U(1\otimes 1\otimes C)$ such that 
\begin{equation}\label{Z}
\text{$\|x-\text{E}_{1\otimes B\otimes C}(x)\|_2\leq \eta\big(\|[Z_1,x]\|_2+\|[Z_2,x]\|_2\big)$, for every $x\in M$,}
\end{equation}
and
\begin{equation}\label{T}
\text{$\|x-\text{E}_{A\otimes B\otimes 1}(x)\|_2\leq\eta\big(\|[T_1,x]\|_2+\|[T_2,x]\|_2\big)$, for every $x\in M$.}
\end{equation}
Let $Z_3\in\mathcal U(A\otimes B\otimes 1)$ and $T_3\in\mathcal U(1\otimes B\otimes C)$ be given by 
$$\text{$Z_3=\sum_{i=1}^ke_{i,i}\otimes U_i\otimes 1$\;\; and \;\; $T_3=\sum_{j=1}^m 1\otimes V_j\otimes e_{j,j}$.}$$
Then $[Z_3,T_3]=\sum_{i=1}^k\sum_{j=1}^me_{i,i}\otimes [U_i,V_j]\otimes e_{j,j}$ and thus $\|[Z_3,T_3]\|_2^2=\frac{1}{km}\sum_{i=1}^k\sum_{j=1}^m\|[U_i,V_j]\|_2^2\leq\delta$. On the other hand, $[Z_{\alpha},T_{\beta}]=0$ if $\alpha,\beta\in\{1,2,3\}$ are not both equal to $3$. Altogether, we get that \begin{equation}\label{F3F3}\text{$\|[Z_{\alpha},T_{\beta}]\|_2^2\leq\delta$,\;\; for every $\alpha,\beta\in\{1,2,3\}$.} \end{equation}

The second paragraph of the proof implies that there are $\widetilde Z_\alpha,\widetilde T_\beta\in\mathcal U(M)$ such that  $\|\widetilde Z_\alpha-Z_\alpha\|_2\leq\varepsilon_0, \|\widetilde T_\beta-T_\beta\|_2\leq\varepsilon_0$ and $[\widetilde Z_\alpha,\widetilde T_\beta]=0$, for all $\alpha,\beta\in\{1,2,3\}$. 

Denote by $P\subset M$ the von Neumann subalgebra generated by $\{\widetilde Z_1,\widetilde Z_2,\widetilde Z_3\}$. 
Let $x\in (P)_1$. If $\beta\in\{1,2\}$, then since $x$ commutes with $\widetilde T_\beta$ we get that $\|[T_\beta,x]\|_2\leq 2\|\widetilde T_\beta-T_\beta\|_2\leq 2\varepsilon_0$ and \eqref{T} gives that $\|x-\text{E}_{A\otimes B\otimes 1}(x)\|_2\leq 4\eta\varepsilon_0$. As $x\in (P)_1$ is arbitrary, we get $P\subset_{4\eta\varepsilon_0}A\otimes B\otimes 1$. 

Similarly, using that $Q=P'\cap M$ commutes with $\widetilde Z_1,\widetilde Z_2$ and \eqref{T} we get that $Q\subset_{4\eta\varepsilon_0}1\otimes B\otimes C$. Since $M$ is a finite dimensional factor, the bicommutant theorem gives that $Q'\cap M=P$. By applying Lemma \ref{commutant} we derive that $A\otimes 1\otimes 1\subset_{16\eta\varepsilon_0}P$.

Since $16\eta\varepsilon_0<\delta_3(\frac{\varepsilon}{16})$,
by combining the last two paragraphs and Lemma \ref{perturb} we find a von Neumann subalgebra $S\subset B$ such that 
\begin{equation}
\label{S}\text{\bf{d}}(P,A\otimes S\otimes 1)\leq \frac{\varepsilon}{16}.
\end{equation}
Denote $T=S'\cap B$. By Lemma \ref{commutant} we get that
\begin{equation}
\label{rel}\text{\bf{d}}(Q,1\otimes T\otimes C)\leq \frac{\varepsilon}{4}.
\end{equation}
Since $\widetilde Z_3\in P$, \eqref{S} gives that $\|\widetilde Z_3-\text{E}_{A\otimes S\otimes 1}(\widetilde Z_3)\|_2\leq \frac{\varepsilon}{16}$.
As $\|Z_3-\widetilde Z_3\|_2\leq\varepsilon_0$, we get that $\|Z_3-\text{E}_{A\otimes S\otimes 1}(Z_3)\|_2\leq \frac{\varepsilon}{16}+2\varepsilon_0<\frac{\varepsilon}{3}$. Similarly, by \eqref{rel} we get $\|T_3-\text{E}_{1\otimes T\otimes C}(T_3)\|_2\leq\frac{\varepsilon}{4}+2\varepsilon_0<\frac{\varepsilon}{3}$. 
The last two inequalities imply that $$\text{$\frac{1}{k}\sum_{i=1}^k\|U_i-\text{E}_S(U_i)\|_2^2\leq\frac{\varepsilon}{9}$\;\; and \;\; $\frac{1}{m}\sum_{j=1}^m\|V_j-\text{E}_T(V_j)\|_2^2\leq\frac{\varepsilon}{9}$.}$$

Finally, by Lemma \ref{close} we can find $\widetilde U_i\in\mathcal U(S),\widetilde V_j\in\mathcal U(T)$ such that $\|U_i-\widetilde U_i\|_2\leq 3\|U_i-\text{E}_S(U_i)\|_2$ and $\|V_j-\widetilde V_j\|_2\leq 3\|V_j-\text{E}_T(V_j)\|_2$, for every $1\leq i\leq k$ and $1\leq j\leq m$. Since $S$ and $T$ commute, the conclusion follows. \hfill$\blacksquare$

\subsection*{Proof of Proposition \ref{reduction}} Assume that $\mathbb F_2\times\mathbb F_2$ is \text{HS} stable. Let $\varepsilon>0$. 
Let $\eta>0$ be the constant provided by Proposition \ref{spF2}.
Let $\delta_3:(0,+\infty)\rightarrow (0,+\infty)$ be the function provided by Lemma \ref{perturb}. Let $\varepsilon_0>0$ such that $\varepsilon_0<\frac{\varepsilon}{24}$ and $16\eta\varepsilon_0<\delta_3(\frac{\varepsilon}{32})$.

Since $\mathbb F_2\times\mathbb F_2$ is \text{HS} stable we can find $\delta>0$ such that for any finite dimensional factor $M$ and $Z_\alpha,T_\beta\in\mathcal U(M)$ such that $\|[Z_\alpha,T_\beta]\|_2^2\leq\delta$, for every $\alpha,\beta\in\{1,2\}$, we can find $\widetilde Z_\alpha,\widetilde T_\beta\in\mathcal U(M)$ such that $\|\widetilde Z_\alpha-Z_\alpha\|_2\leq\varepsilon, \|\widetilde T_\beta-T_\beta\|_2\leq\varepsilon$ and $[\widetilde Z_\alpha,\widetilde T_\beta]=0$, for every $\alpha,\beta\in\{1,2\}$.

Let $k,m,n\in\mathbb N$ and $U_1,...,U_k,V_1,...,V_m\in U(n)$ such that \begin{equation}\label{almostc}\frac{1}{km}\sum_{i=1}^k\sum_{j=1}^m\|[U_i,V_j]\|_2^2\leq\delta.\end{equation} 
Denote $A=\mathbb M_k(\mathbb C), B=\mathbb M_n(\mathbb C)$, $C=\mathbb M_m(\mathbb C)$ and $M=\mathbb M_3(\mathbb C)\otimes A\otimes B\otimes C\otimes\mathbb M_3(\mathbb C).$
Let $W\in\mathcal U(A\otimes B)$ and $W'\in\mathcal U(B\otimes C)$ be given by $W=\sum_{i=1}^ke_{i,i}\otimes U_i$ and $W'=\sum_{j=1}^mV_j\otimes e_{j,j}$.
By applying Lemma \ref{spF2} twice we can find $Z_1,Z_2,T_1,T_2\in\mathcal U(M)$ such that
\begin{enumerate}
\item\label{Z1} $\|x-\text{E}_{1\otimes 1\otimes B\otimes C\otimes\mathbb M_3(\mathbb C)}(x)\|_2\leq\eta(\|[Z_1,x]\|_2+\|[Z_2,x]\|_2)$, for every $x\in M$,
\item $Z_1\in\mathbb M_3(\mathbb C)\otimes A\otimes 1\otimes 1\otimes 1$,  
\item $Z_2=\begin{pmatrix} 0&0&1\\1&0&0\\0&W&0 \end{pmatrix}\otimes 1\otimes 1\in \mathbb M_3(\mathbb C)\otimes A\otimes B\otimes 1\otimes 1$, 
\item\label{T1} $\|x-\text{E}_{\mathbb M_3(\mathbb C)\otimes A\otimes B\otimes 1\otimes 1}(x)\|_2\leq\eta(\|[T_1,x]\|_2+\|[T_2,x]\|_2)$, for every $x\in M$,
\item $T_1\in 1\otimes 1\otimes 1\otimes C\otimes\mathbb M_3(\mathbb C)$ and
\item $T_2=1\otimes 1\otimes \begin{pmatrix} 0&0&1\\1&0&0\\0&W'&0 \end{pmatrix}\in 1\otimes 1\otimes B\otimes C\otimes\mathbb M_3(\mathbb C)$.

\end{enumerate}

Next, we have $[Z_2,T_2]=e_{3,2}\otimes [W\otimes 1,1\otimes W']\otimes e_{3,2}=\frac{1}{km}\sum_{i=1}^k\sum_{j=1}^me_{3,2}\otimes e_{i,i}\otimes [U_i,V_j]\otimes e_{j,j}\otimes e_{3,2}$
and thus $\|[Z_2,T_2]\|_2^2=\frac{1}{9km}\sum_{i=1}^k\sum_{j=1}^m\|[U_i,V_j]\|_2^2\leq \delta.$ On the other hand, $[Z_\alpha,T_\beta]=0$ if $\alpha,\beta\in\{1,2\}$ are not both equal to $2$. Altogether, we have that  $\|[Z_\alpha,T_\beta]\|_2^2\leq\delta$, for every $\alpha,\beta\in\{1,2\}$. 

The second paragraph of the proof implies that we can find $\widetilde Z_\alpha,\widetilde T_\beta\in\mathcal U(M)$ such that $\|\widetilde Z_\alpha-Z_\alpha\|_2\leq\varepsilon_0, \|\widetilde T_\beta-T_\beta\|_2\leq\varepsilon_0$ and $[\widetilde Z_\alpha,\widetilde T_\beta]=0$, for every $\alpha,\beta\in\{1,2\}$. 

Let $P\subset M$ be the von Neumann subalgebra generated by $\{Z_1,Z_2\}$. Let $x\in (P)_1$. If $\beta\in\{1,2\}$, then since $x$ commutes with $\widetilde T_\beta$, we get that $\|[T_\beta,x]\|_2\leq 2\|\widetilde T_\beta-T_\beta\|_2\leq 2\varepsilon_0$ and \eqref{T1} gives that $\|x-\text{E}_{\mathbb M_3(\mathbb C)\otimes A\otimes B\otimes 1\otimes 1}(x)\|_2\leq 4\eta\varepsilon_0$. 
As $x\in (P)_1$ is arbitrary, we get $P\subset_{4\eta\varepsilon_0}\mathbb M_3(\mathbb C)\otimes A\otimes B\otimes 1\otimes 1$.

Similarly, using that $Q=P'\cap M$ commutes with $\widetilde Z_1,\widetilde Z_2$ and \eqref{Z1}, we deduce that $Q\subset_{4\eta\varepsilon_0}1\otimes 1\otimes B\otimes C\otimes\mathbb M_3(\mathbb C)$.
Since $M$ is a finite dimensional factor, the bicommutant theorem gives that $Q'\cap M=P$. By applying Lemma \ref{commutant} we get that $\mathbb M_3(\mathbb C)\otimes A\otimes 1\otimes 1\otimes 1\subset_{16\eta\varepsilon_0}P$. 

Since $16\eta\varepsilon_0<\delta_3(\frac{\varepsilon}{32})$, using the last two paragraphs and Lemma \ref{perturb} we find a von Neumann subalgebra $S\subset B$ such that 
\begin{equation}
\label{S2}\text{\bf{d}}(P,\mathbb M_3(\mathbb C)\otimes A\otimes S\otimes 1\otimes 1)\leq \frac{\varepsilon}{32}.
\end{equation}
Denote $T=S'\cap B$. By Lemma \ref{commutant} we get that
\begin{equation}
\label{rel2}\text{\bf{d}}(Q,1\otimes 1\otimes  T\otimes C\otimes \mathbb M_3(\mathbb C))\leq \frac{\varepsilon}{8}.
\end{equation}

Since $\widetilde Z_2\in P$, \eqref{S2} gives that $\|\widetilde Z_2-\text{E}_{\mathbb M_3(\mathbb C)\otimes A\otimes S\otimes 1\otimes 1}(\widetilde Z_2)\|_2\leq \frac{\varepsilon}{32}$.
Since $\|Z_2-\widetilde Z_2\|_2\leq\varepsilon_0$, we get that $\|Z_2-\text{E}_{\mathbb M_3(\mathbb C)\otimes A\otimes S\otimes 1\otimes 1}(Z_2)\|_2\leq \frac{\varepsilon}{32}+2\varepsilon_0<\frac{\varepsilon}{6}$. Similarly, using \eqref{rel} we get that $\|T_2-\text{E}_{1\otimes 1\otimes T\otimes C\otimes\mathbb M_3(\mathbb C)}(T_2)\|_2\leq \frac{\varepsilon}{8}+2\varepsilon_0<\frac{\varepsilon}{6}$. 
By the definition of $Z_2, T_2$, the last two inequalities imply that $$\text{$\frac{1}{3k}\sum_{i=1}^k\|U_i-\text{E}_S(U_i)\|_2^2\leq \frac{\varepsilon}{36}$ and $\frac{1}{3m}\sum_{j=1}^m\|V_j-\text{E}_T(V_j)\|_2^2\leq \frac{\varepsilon}{36}$.}$$

Finally, by Lemma \ref{close} we can find $\widetilde U_i\in\mathcal U(S),\widetilde V_j\in\mathcal U(T)$ such that $\|U_i-\widetilde U_i\|_2\leq 3\|U_i-\text{E}_S(U_i)\|_2$ and $\|V_j-\widetilde V_j\|_2\leq 3\|V_j-\text{E}_T(V_j)\|_2$, for every $1\leq i\leq k$ and $1\leq j\leq m$. Since $S$ and $T$ commute, the conclusion follows. \hfill$\blacksquare$

\section{Proof of Theorem \ref{main}}

\subsection{Construction}
In this section we prove Theorem \ref{main} by constructing a counterexample to the conclusion of Proposition \ref{reduction}.   We start by recalling our construction presented in the introduction.

\begin{notation}\label{construction}
Let $n\in\mathbb N$ and $t\in\mathbb R$. 
\begin{enumerate}
\item We denote $M_n=\bigotimes_{k=1}^n\mathbb M_2(\mathbb C)\cong\mathbb M_{2^n}(\mathbb C)$ and $A_n=\bigotimes_{k=1}^n\mathbb C^2\cong \mathbb C^{2^n}$. We view $A_n$ as a subalgebra of $M_n$, where we embed $\mathbb C^2\subset\mathbb M_2(\mathbb C)$ as the diagonal matrices.

 \item For $1\leq i\leq n$, let $X_{n,i}=1\otimes\cdots 1\otimes \sigma\otimes 1 \cdots\otimes 1\in \mathcal U(A_n)$, where $\sigma=\begin{pmatrix} 1&0\\ 0&-1\end{pmatrix}\in\mathbb C^2$ is placed on the $i$-th tensor position.

 \item Let $G_n\subset\mathcal U(A_n\otimes M_n)$ be a finite subgroup which generates $A_n\otimes M_n$. 


 \item We define $U_t\in\text{U}(\mathbb C^2\otimes\mathbb C^2)$ by  $U_t=P+e^{it}(1-P)$, where $P:\mathbb C^2\otimes\mathbb C^2\rightarrow\mathbb C^2\otimes\mathbb C^2$ be the orthogonal projection onto the one dimensional space spanned by $e_1\otimes e_2-e_2\otimes e_1$ .

 \item We identify $M_n\otimes M_n=\otimes_{k=1}^n(\mathbb M_2(\mathbb C)\otimes\mathbb M_2(\mathbb C))$ and  let $\theta_{t,n}$ be the automorphism of $M_n\otimes M_n$ given by $\theta_{t,n}(\otimes_{k=1}^nx_k)=\otimes_{k=1}^nU_tx_kU_t^*$. 
\item Finally, consider the following two sets of unitaries in $M_n\otimes M_n$: $\mathcal U_n=\{X_{n,i}\otimes 1\mid 1\leq i\leq n\}$ and $\mathcal V_{t,n}=G_n\cup\theta_{t,n}(G_n)$.
\end{enumerate}
\end{notation}

We begin with the following elementary lemma. For $t\in \mathbb R$, we let $\rho_t=\frac{1+\cos(t)}{2}\in [0,1]$.  We endow $M_n$ with its unique trace $\tau$ and the scalar product given by $\langle x,y\rangle=\tau(y^*x)$, for every $x,y\in M_n$. For $1\leq l\leq n$, we denote by  $e_l:M_n\rightarrow M_n$ the orthogonal projection onto the subspace of tensors of length at most $l$, i.e., the span of $\otimes_{k=1}^nx_k$, with $x_k\in\mathbb M_2(\mathbb C)$ and $|\{k\mid x_k\not=1\}|\leq l$.

\begin{lemma}\label{deform} 
The following hold:
\begin{enumerate}
\item $P(x\otimes 1)P=\tau(x)P$, for every $x\in\mathbb M_2(\mathbb C)$.
\item $\tau((x\otimes 1)U_t(y\otimes 1)U_t^*)=\rho_t\tau(xy)+(1-\rho_t)\tau(x)\tau(y)$, for every $x,y\in\mathbb M_2(\mathbb C)$.
\item $\emph{E}_{U_t(\mathbb M_2(\mathbb C)\otimes 1)U_t^*}(x\otimes 1)=\rho_t U_t(x\otimes 1)U_t^*$, for every $x\in\mathbb M_2(\mathbb C)$ with $\tau(x)=0$.
\item $\|\emph{E}_{\theta_{t,n}(M_n\otimes 1)}(x\otimes 1)\|_2^2\leq (1-\rho_t^{2l}) \|e_l(x)\|_2^2+\rho_t^{2l}\|x\|_2^2$, for every $x\in M_n$ and $1\leq l\leq n$.
\end{enumerate}
\end{lemma}

{\it Proof.} It is immediate that $P(e_{i,j}\otimes 1)P$ is equal to $\frac{1}{2}P$, if $i=j$, and $0$, if $i\not=j$, which implies (1). 
Part (2) follows via a straightforward calculation by using (1) and that $\tau(P)=\frac{1}{4}$. 
If $x\in\mathbb M_2(\mathbb C)$ and $\tau(x)=0$, then (2) gives $\tau((x\otimes 1)U_t(y\otimes 1)U_t^*)=\rho_t\tau(xy)=\rho_t\tau(U_t(x\otimes 1)U_t^*U_t(y\otimes 1)U_t^*),$ for every $y\in\mathbb M_2(\mathbb C)$. This clearly implies (3).

To prove (4), for $0\leq i\leq n$, we denote by $V_i\subset M_n$ the span of tensors of the form $\otimes _{k=1}^nx_k$,  such that $x_k=1$ or $\tau(x_k)=0$, for every $1\leq k\leq n$, and $|\{k\mid x_k\not=1\}|=i$. Let $f_i:M_n\rightarrow M_n$ be the orthogonal projection onto $V_i$. If $x=\otimes_{k=1}^nx_k\in V_i$, then using part (3) we get that  $$\text{E}_{\theta_{t,n}(M_n\otimes 1)}(x\otimes 1)=\otimes_{k=1}^n\text{E}_{U_t(\mathbb M_2(\mathbb C)\otimes 1)U_t^*}(x_k\otimes 1)=\rho_t^{i}\theta_{t,n}(x\otimes 1).$$
Thus, for every $x\in M_n$ we have $\text{E}_{\theta_{t,n}(M_n\otimes 1)}(x\otimes 1)=\sum_{i=0}^n\rho_t^i\theta_{t,n}(f_i(x)\otimes 1)$ and therefore
\begin{equation}\label{theta_n}\|\text{E}_{\theta_{t,n}(M_n\otimes 1)}(x\otimes 1)\|_2^2=\sum_{i=0}^n\rho_t^{2i}\|f_i(x)\|_2^2.\end{equation}
Since  $\sum_{i=0}^l\|f_i(x)\|_2^2=\|e_l(x)\|_2^2$ and $\sum_{i=l+1}^n\|f_i(x)\|_2^2=\|x|_2^2-\|e_l(x)\|_2^2$, \eqref{theta_n} implies part (4). 
\hfill$\blacksquare$

Next, we show that the sets of unitaries $\mathcal U_n$ and $\mathcal V_{t,n}$ almost commute:
\begin{lemma}\label{almostcommute}
$\|[U,V]\|_2\leq \|[U,V]\|\leq 4|t|$, for all $U\in\mathcal U_n$, $V\in\mathcal V_{t,n}$.
\end{lemma}

{\it Proof.}
Let $U\in\mathcal U_n$ and $V\in\mathcal V_{t,n}$. If $V\in G_n\subset A_n\otimes M_n$, then as $U\in A_n\otimes 1$ and $A_n$ is abelian we get $[U,V]=0$.
Thus, we may assume that $V=\theta_{t,n}(Y)$, for $Y\in G_n$. Then since $[\theta_{t,n}(U),V]=0$, we get  $\|[U,V]\|\leq 2\|U-\theta_{t,n}(U)\|=2\|(\sigma\otimes 1)-U_t(\sigma\otimes 1)U_t^*\|\leq 4\|U_t-1\|=4|e^{it}-1|\leq 4|t|$.
\hfill$\blacksquare$

\subsection{A consequence of HS-stability of $\mathbb F_2\times\mathbb F_2$}
To prove Theorem \ref{main}, we show that if $t>0$ is small enough, the almost commuting sets of unitaries $\mathcal U_n$ and $\mathcal V_{t,n}$ contradict the conclusion of Proposition \ref{reduction} for large $n\in\mathbb N$.
To this end, we first use Proposition \ref{reduction} to deduce the following:
\begin{corollary}\label{coro}
Assume that $\mathbb F_2\times\mathbb F_2$ is \emph{HS}-stable. Then for every $\varepsilon>0$, there exists $t>0$ such that the following holds: for every $n\in\mathbb N$, we can find a von Neumann subalgebra $C\subset A_n$ such that 
\begin{enumerate}
\item $C\otimes 1\subset_{\varepsilon}\theta_{t,n}(A_n\otimes 1)$ and
\item $\frac{1}{n}\sum_{i=1}^n\|X_{n,i}-\emph{E}_C(X_{n,i})\|_2^2\leq\varepsilon$.
\end{enumerate}
\end{corollary}

{\it Proof.}
Let $\varepsilon\in (0,1)$. Let $\eta>0$ such that $\eta<\frac{\varepsilon^2}{256}$ and $\eta<\frac{\delta_2(\frac{\varepsilon}{2})^2}{64}$, where $\delta_2:(0,+\infty)\rightarrow (0,+\infty)$ is the function provided by Corollary \ref{abelian}. 

 By Lemma  \ref{almostcommute} we have that $\frac{1}{|\mathcal U_n|\cdot |\mathcal V_{t,n}|}\sum_{U\in\mathcal U_n,V\in\mathcal V_{t,n}}\|[U,V]\|_2^2\leq 16t^2$, for every $n\in\mathbb N$ and $t\in\mathbb R$. 
Since $\mathbb F_2\times\mathbb F_2$ is HS-stable, Proposition \ref{reduction} implies that if $t>0$ is small enough then the following holds: given any $n\in\mathbb N$,
 we can find a von Neumann subalgebra $P\subset M_n\otimes M_n$ such that 
 \begin{equation}\label{Un}\text{$\frac{1}{|\mathcal U_n|}\sum_{U\in\mathcal U_n}\|U-\text{E}_P(U)\|_2^2\leq\eta$}\end{equation} and
\begin{equation}\label{Vtn} \frac{1}{|\mathcal V_{t,n}|}\sum_{V\in\mathcal V_{t,n}}\|V-\text{E}_{P'}(V)\|_2^2\leq\eta.
\end{equation}
Then \eqref{Vtn} gives that $\frac{1}{|G_n|}\sum_{V\in G_n}\|V-\text{E}_{P'}(V)\|_2^2\leq 2\eta$ and $\frac{1}{|G_n|}\sum_{V\in \theta_{t,n}(G_n)}\|V-\text{E}_{P'}(V)\|_2^2\leq 2\eta$.
Since $G_n$ generates $A_n\otimes M_n$, by Lemma \ref{basic} we conclude that  \begin{equation}\label{AnMn}\text{$A_n\otimes M_n\subset_{2\sqrt{\eta}}P'$\;\;\; and\;\;\; $\theta_{t,n}(A_n\otimes M_n)\subset_{2\sqrt{\eta}}P'$.}\end{equation}
Since $M_n\otimes M_n$ is a finite dimensional factor, the bicommutant theorem gives that $(P')'=P$. Since $A_n\subset M_n$ is a maximal abelian subalgebra, we also have that $(A_n\otimes M_n)'=A_n\otimes 1$. By combining these facts with \eqref{AnMn} and Lemma \ref{commutant}, we derive that $P\subset_{8\sqrt{\eta}}A_n\otimes 1$ and $P\subset_{8\sqrt{\eta}}\theta_{t,n}(A_n\otimes 1)$. 

Since $A_n\otimes 1$ is abelian and we have chosen $\eta>0$ so that $8\sqrt{\eta}\leq\delta_2(\frac{\varepsilon}{2})$, Corollary \ref{abelian} implies that we can find a von Neumann subalgebra $Q\subset A_n\otimes 1$ such that $\text{\bf{d}}(P,Q)\leq\frac{\varepsilon}{2}$. Since $8\sqrt{\eta}\leq \frac{\varepsilon}{2}$, we also have that $P\subset_{\frac{\varepsilon}{2}}\theta_{t,n}(A_n\otimes 1)$. By combining the last two facts we derive that $Q\subset_{\varepsilon}\theta_{t,n}(A_n\otimes 1)$.

Thus, if $C\subset A_n$ is a von Neumann subalgebra such that $Q=C\otimes 1$, then condition (1) is satisfied. 
To verify condition (2), let $U\in\mathcal U(M_n\otimes M_n)$.
Since $\text{\bf{d}}(P,Q)\leq\frac{\varepsilon}{2}$ we have that \begin{align*}\|U-\text{E}_Q(U)\|_2&\leq \|U-\text{E}_P(U)\|_2+\|\text{E}_P(U)-\text{E}_Q(\text{E}_P(U))\|_2+\|\text{E}_Q(\text{E}_P(U)-U)\|_2\\&\leq 2\|U-\text{E}_P(U)\|_2+\frac{\varepsilon}{2}.
\end{align*}
Hence, $\|U-\text{E}_Q(U)\|_2^2\leq 2(4\|U-\text{E}_P(U)\|_2^2+\frac{\varepsilon^2}{4})=8\|U-\text{E}_P(U)\|_2^2+\frac{\varepsilon^2}{2}$. In combination with \eqref{Un}, we derive that $\frac{1}{n}\sum_{i=1}^n\|X_{n,i}-\emph{E}_C(X_{n,i})\|_2^2=\frac{1}{|\mathcal U_n|}\sum_{U\in\mathcal U_n}\|U-\text{E}_Q(U)\|_2^2\leq 8\eta+\frac{\varepsilon^2}{2}$. Since $\eta<\frac{\varepsilon^2}{256}$ and $\varepsilon\in (0,1)$ we have that $8\eta+\frac{\varepsilon^2}{2}<\varepsilon$ and  condition (2) follows. 
\hfill$\blacksquare$

Let $\varepsilon\in(0,\frac{1}{16})$. Assuming that $\mathbb F_2\times\mathbb F_2$ is HS-stable, Corollary \ref{coro} implies that there is $t>0$ such that for every $n\in\mathbb N$ we can find a subalgebra  $C_n\subset A_n$ such that (a)  $C_n\otimes 1\subset_{\varepsilon}\theta_{t,n}(A_n\otimes 1)$ and (b) $\frac{1}{n}\sum_{i=1}\|X_{n,i}-\text{E}_{C_n}(X_{n,i})\|_2^2<\varepsilon$. We will derive a contradiction as $n\rightarrow\infty$ by showing that (a) and (b) imply the following incompatible facts:
\begin{itemize}
\item $\text{dim}(C_nz_n)\leq P(n)$, where $P$ is a polynomial (see Lemma \ref{1}), and
\item $\text{dim}(C_nz_n)\geq 2^{\kappa n}$, where $\kappa\in (0,1)$ (see Lemma \ref{2}), for a projection $z_n\in C_n$.
\end{itemize}

\subsection{A polynomial upper bound on dimension}

\begin{lemma}\label{1}
Let $C\subset M_n$ be a von Neumann subalgebra such that $C\otimes 1\subset_{\varepsilon}\theta_{t,n}(M_n\otimes 1)$, for some $n\in\mathbb N$,  $t\in (0,\frac{\pi}{4}]$ and $\varepsilon\in [0,\frac{1}{16}]$.

Then there exists a projection $z\in\mathcal Z(C)$ such that $\emph{dim}(Cz)\leq 2(6n)^{64\frac{\varepsilon}{t^2}+1}$ and $\tau(z)\geq\frac{1}{2}$.
\end{lemma}

{\it Proof.} 
For simplicity, we denote $D=\theta_{t,n}(M_n\otimes 1)$. Since $C\otimes 1\subset_{\varepsilon}D$, we get that
\begin{equation}\label{E_D}
\text{$\|\text{E}_D(u\otimes 1)\|_2\geq 1-\varepsilon$, for every $u\in\mathcal U(C)$.}
\end{equation}

Let $\{z_j\}_{j=1}^m$ be an enumeration of the minimal projections of $\mathcal Z(C)$. Then $C=\bigoplus_{j=1}^mCz_j$, where $Cz_j$ is a factor and thus isomorphic to a matrix algebra $\mathbb M_{n_j}(\mathbb C)$, for some $n_j\in\mathbb N$.

Let $S$ be the set of $1\leq j\leq m$ such that $\|E_D(u\otimes 1)\|_2^2\geq (1-4\varepsilon)\|u\|_2^2=(1-4\varepsilon)\tau(z_j)$, for every $u\in\mathcal U(Cz_j)$.
Let $T=\{1,\cdots,m\}\setminus S$. Then for every $j\in T$,  there exists $u_j\in\mathcal U(Cz_j)$ such that \begin{equation}\label{u_j}\|\text{E}_D(u_j\otimes 1)\|_2^2\leq (1-4\varepsilon)\tau(z_j).\end{equation}
We will prove that $z=\sum_{j\in S}z_j\in\mathcal Z(C)$ satisfies the conclusion. To estimate $\tau(z)$, for every $j\in S$, let $u_j=1$.
Denote by $\lambda^m$ the Haar measure of $\mathbb T^m$, where $\mathbb T=\{z\in\mathbb C\mid |z|=1\}$. By applying \eqref{E_D} to $\sum_{j=1}^m\mu_ju_j\in\mathcal U(C)$, for $\mu_1,...,\mu_m\in\mathbb T$, we get that 
\begin{equation}\label{integration}
\sum_{j=1}^m\|\text{E}_D(u_j\otimes 1)\|_2^2=\int_{\mathbb T^m}\|\text{E}_D(\sum_{j=1}^m\mu_ju_j\otimes 1)\|_2^2\;\text{d}\lambda^m(\mu_1,\cdots,\mu_m)\geq (1-\varepsilon)^2\geq 1-2\varepsilon.
\end{equation}
On the other, \eqref{u_j} implies that $\sum_{j=1}^m\|\text{E}_D(u_j\otimes 1)\|_2^2\leq \tau(z)+(1-4\varepsilon)\tau(1-z)$. In combination with \eqref{integration} we deduce that $\tau(z)\geq\frac{1}{2}$.

To estimate $\text{dim}(Cz)$, let $l$ be the smallest positive integer such that $\rho_t^{2l}\leq 1-8\varepsilon$.
We claim that
\begin{equation}\label{1/2}
\text{$\|e_l(u)\|_2^2\geq\frac{\|u\|_2^2}{2}$, for every $j\in S$ and $u\in\mathcal U(Cz_j)$.}
\end{equation}

If $u\in\mathcal U(Cz_j)$, for some $j\in S$, then  Lemma \ref{deform}(4) gives that
\begin{align*}(1-4\varepsilon)\|u\|_2^2&\leq \|\text{E}_D(u\otimes 1)\|_2^2\\&\leq \rho_t^{2l}(\|u\|_2^2-\|e_l(u)\|_2^2)+\|e_l(u)\|_2^2\\&\leq (1-8\varepsilon)(\|u\|_2^2-\|e_l(u)\|_2^2)+\|e_l(u)\|_2^2,\end{align*} which implies \eqref{1/2}.

If $j\in S$, then since $Cz_j$ is a isomorphic to the matrix algebra $\mathbb M_{n_j}(\mathbb C)$, it admits an orthonormal basis $\mathcal B_j$ whose every element is of the form $\frac{u}{\|u\|_2}$, for some $u\in\mathcal U(Cz_j)$. Then $\mathcal B=\cup_{j\in S}\mathcal B_j$  is an orthonormal basis for $Cz=\bigoplus_{j\in S}Cz_j$ and \eqref{1/2} implies that $\|e_l(\xi)\|_2^2\geq \frac{1}{2}$, for every $\xi\in\mathcal B$. Recall that $e_l$ is the orthogonal projection onto the subspace $W_l\subset M_n$ of tensors of length at most $l$ and let $\mathcal O$ be an orthonormal basis for $W_l$. Then we have that \begin{equation}\label{dim}
\text{dim}(Cz)=|\mathcal B|\leq 2\sum_{\xi\in\mathcal B}\|e_l(\xi)\|_2^2=2\sum_{\xi\in\mathcal B,\eta\in\mathcal O}|\langle\xi,\eta\rangle|^2\leq 2|\mathcal O|=2\;\text{dim}(W_l).
\end{equation}
On the other hand, we have the following crude estimate
\begin{equation}
\label{dimm}
\text{dim}(W_l)=\sum_{i=0}^l3^i{n\choose i}\leq (l+1)3^ln^l\leq (6n)^l.
\end{equation}
Next, note that $x\leq |\log(1-x)|\leq 2x$, for every $x\in [0,\frac{1}{2}]$. Since $\varepsilon\in [0,\frac{1}{16}]$, we get that $|\log(1-8\varepsilon)|\leq 16\varepsilon$. Since $t\in (0,\frac{\pi}{4}]$, we also have that $1-\rho_t=\frac{1-\cos(t)}{2}\in [0,\frac{1}{2}]$ and $1-\rho_t\geq\frac{t^2}{8}$. Thus, $|\log(\rho_t)|\geq 1-\rho_t\geq\frac{t^2}{8}$.
By using these facts and the definition of $l$ we derive that \begin{equation}\label{ll}l\leq\frac{|\log(1-8\varepsilon)|}{2|\log(\rho_t)|}+1\leq 64\frac{\varepsilon}{t^2}+1.\end{equation}
Combining \eqref{dim}, \eqref{dimm} and \eqref{ll} implies that $\text{dim}(Cz)\leq 2(6n)^{64\frac{\varepsilon}{t^2}+1}$, as desired. \hfill$\blacksquare$

\subsection{An exponential lower bound on dimension}
Let $\text{H}:(0,1)\rightarrow (0,1]$ be the binary entropy function given by $\text{H}(\delta)=-\delta\log_2(\delta)-(1-\delta)\log_2(1-\delta)$.

\begin{lemma}\label{2}
Let $C\subset A_n$ be a von Neumann subalgebra such that $\frac{1}{n}\sum_{i=1}^n\|X_{n,i}-\emph{E}_C(X_{n,i})\|_2^2\leq \varepsilon$, for some $\varepsilon\in [0,\frac{1}{8}]$ and $n\in\mathbb N$. 
Then $\emph{dim}(Cz)\geq 2^{n-\emph{H}(4\varepsilon)n-3}$, for any projection $z\in C$ with $\tau(z)\geq\frac{1}{2}$.
 \end{lemma}

{\it Proof.} Let $p=\begin{pmatrix} 1&0\\0&0\end{pmatrix}$. For $1\leq i\leq n$, let $p_{i}=1\otimes \cdots 1\otimes p\otimes 1\cdots 1$, where $p$ is placed on the $i$-th tensor position. Then $X_{n,i}=2p_{i}-1$ and so $X_{n,i}-\text{E}_C(X_{n,i})=2(p_i-\text{E}_C(p_i))$, for every $1\leq i\leq n$. Thus, the hypothesis rewrites as 
\begin{equation}\label{contain}\sum_{i=1}^n\|p_i-\text{E}_C(p_i)\|_2^2\leq\frac{\varepsilon}{4}.
\end{equation}

Let $\{q_j\}_{j=1}^m$ be the minimal projections of $C$ such that $C=\bigoplus_{j=1}^m\mathbb Cq_j$. We claim that 
\begin{equation}\label{expect}
\text{$\|p-\text{E}_C(p)\|_2^2=\sum_{j=1}^m\frac{\tau(pq_j)\tau((1-p)q_j)}{\tau(q_j)}$, for every projection $p\in A_n$.}\end{equation}

Since $\text{E}_C(p)=\sum_{j=1}^m\frac{\tau(pq_j)}{\tau(q_j)}q_j$, we get that $\|\text{E}_C(p)\|_2^2=\sum_{j=1}^m\frac{\tau(pq_j)^2}{\tau(q_j)}$. By combining the last fact with the identity $\|p-E_C(p)\|_2^2=\|p\|_2^2-\|E_C(p)\|_2^2=\tau(p)-\|E_C(p)\|_2^2$, \eqref{expect} follows.

By combining \eqref{contain} and \eqref{expect} we deduce that
\begin{equation}\label{piqj}\frac{1}{n}\sum_{i=1}^n\sum_{j=1}^m\frac{\tau(p_iq_j)\tau((1-p_i)q_j)}{\tau(q_j)}\leq\frac{\varepsilon}{4}.
\end{equation}

Next, let $S$ be the set of $1\leq j\leq m$ such that $\frac{1}{n}\sum_{i=1}^n\tau(p_iq_j)\tau((1-p_i)q_j)<\varepsilon\tau(q_j)^2$.
Define $T=\{1,\cdots,m\}\setminus S$. Let $r=\sum_{j\in T}q_j$. Since $\varepsilon\tau(q_j)\leq \frac{1}{n}\sum_{i=1}^n\frac{\tau(p_iq_j)\tau((1-p_i)q_j)}{\tau(q_j)}$, for every $j\in T$, by using \eqref{piqj} we derive that $\varepsilon\tau(r)=\varepsilon\sum_{j\in T}\tau(q_j)\leq\frac{\varepsilon}{4}$ and thus $\tau(r)\leq\frac{1}{4}$.

\begin{claim}\label{qj}
$\tau(q_j)\leq 2^{H(4\varepsilon)n-n+1}$, for every $j\in S$. 
\end{claim}

{\it Proof of Claim \ref{qj}.}
We identify $A_n$ with $\text{L}^{\infty}(\{0,1\}^n,\mu)$, where $\mu$ is the uniform probability measure on $\{0,1\}^n$. Then $p_i$ is identified with the characteristic function of the set $\{x\in\{0,1\}^n\mid x_i=0\}$, for every $1\leq i\leq n$, and $\tau({\bf 1}_Y)=\mu(Y)=\frac{|Y|}{2^n}$, for every $Y\subset\{0,1\}^n$.

For $x,y\in\{0,1\}$ we denote the normalized Hamming distance:
$$\text{d}_{\text{H}}(x,y)=\frac{|\{i\in\{1,\cdots,n\}\mid x_i\not=y_i\}|}{n}.$$

Let $j\in S$ and $Z_j\subset \{0,1\}^n$ such that $q_j={\bf 1}_{Z_j}$. Since $\tau(p_iq_j)=\mu(\{x\in Z_j\mid x_i=0\})$ and $\tau((1-p_i)q_j)=\mu(\{x\in Z_j\mid x_i=1\})$, the inequality $\varepsilon\tau(q_j)^2>\frac{1}{n}\sum_{i=1}^n\tau(p_iq_j)\tau((1-p_i)q_j)$ rewrites as \begin{align*}\varepsilon\mu(Z_j)^2&>\frac{1}{n}\sum_{i=1}\Big(\mu(\{x\in Z_j\mid x_i=0\})\cdot\mu(\{x\in Z_j\mid x_i=1\})\Big)\\&=\frac{1}{2n}\sum_{i=1}^n(\mu\times\mu)(\{(x,y)\in Z_j\times Z_j\mid x_i\not=y_i\})\\&=\frac{1}{2}\int_{Z_j\times Z_j}\text{d}_{\text{H}}(x,y)\;\text{d}(\mu\times\mu)(x,y).\end{align*}

By Fubini's theorem we can find $x\in Z_j$ such that $\int_{Z_j}\text{d}_{\text{H}}(x,y)\;\text{d}\mu(y)<2\varepsilon\mu(Z_j)$. This implies that $\mu(\{y\in Z_j\mid\text{d}_{\text{H}}(x,y)\geq 4\varepsilon\})<\frac{\mu(Z_j)}{2}$ and hence $\mu(\{y\in Z_j\mid\text{d}_{\text{H}}(x,y)<4\varepsilon\})> \frac{\mu(Z_j)}{2}$.
Thus, 

\begin{align*}
\mu(Z_j)&<2\mu(\{y\in \{0,1\}^n\mid\text{d}_{\text{H}}(x,y)<4\varepsilon\})\\&=2\mu(\{y\in\{0,1\}^n\mid\frac{1}{n}\sum_{i=1}^ny_i<4\varepsilon\})\\&=\frac{1}{2^{n-1}}\sum_{i=0}^{\lfloor 4\varepsilon\rfloor}{n\choose i}.
\end{align*}
Since $\sum_{i=0}^{\lfloor\delta n\rfloor}{n\choose i}\leq 2^{\text{H}(\delta)n}$, for every $n\in\mathbb N$ and $\delta\in (0,\frac{1}{2}]$ (see, e.g., \cite[Lemma 16.19]{FG06}), in combination with the last displayed inequality, we conclude that $\tau(q_j)=\mu(Z_j)\leq 2^{H(4\varepsilon)n-n+1}$.
\hfill$\square$

To finish the proof, let $z\in C$ be a projection with $\tau(z)\geq\frac{1}{2}$. Since $\tau(1-r)=1-\tau(r)\geq\frac{3}{4}$, we have that $\tau(z(1-r))\geq \frac{1}{4}$.  Since $z(1-r)\in C(1-r)=\bigoplus_{j\in S}\mathbb Cq_j$, there is a subset $S_0\subset S$ such that $z(1-r)=\sum_{j\in S_0}q_j$. 
By Claim \eqref{qj} we get that $\frac{1}{4}\leq\tau(z(1-r))=\sum_{j\in S_0}\tau(q_j)\leq |S_0| 2^{H(4\varepsilon)n-n+1}$ and thus $|S_0|\geq 2^{n-H(4\varepsilon)n-3}$.
Since $Cz\supset Cz(1-r)=\bigoplus_{j\in S_0}\mathbb Cq_j$, we have that $\text{dim}(Cz)\geq |S_0|$ and the conclusion follows.
\hfill$\blacksquare$

\subsection{Proof of Theorem \ref{main}} Assume by contradiction that $\mathbb F_2\times\mathbb F_2$ is \text{HS}-stable. Let $\varepsilon\in (0,\frac{1}{16})$.
Then by Corollary \ref{coro} there exists $t>0$ such that for every $n\in\mathbb N$ we can find a von Neumann subalgebra $C_n\subset A_n$ such that (a) $C_n\otimes 1\subset_{\varepsilon}\theta_{t,n}(A_n\otimes 1)$ and  (b) $\frac{1}{n}\sum_{i=1}^n\|X_{n,i}-\text{E}_{C_n}(X_{n,i})\|_2^2\leq\varepsilon$.

Using (a), Lemma \ref{1} gives a projection $z_n\in C_n$ so that $\text{dim}(C_nz_n)\leq 2(6n)^{64\frac{\varepsilon}{t^2}+1}$ and $\tau(z_n)\geq\frac{1}{2}$. On the other hand, using (b), Lemma \ref{2} implies that $\text{dim}(C_nz_n)\geq 2^{n-\text{H}(4\varepsilon)n-3}$. 
Thus, $2(6n)^{64\frac{\varepsilon}{t^2}+1}\geq 2^{n-\text{H}(4\varepsilon)n-3}$, for all $n\in\mathbb N$.
Since $\text{H}(4\varepsilon)<1$, letting $n\rightarrow\infty$ gives a contradiction. \hfill$\blacksquare$

\section{Proofs of Theorem \ref{acm} and Corollary \ref{C}} 
In this section we give the proofs of Theorem \ref{acm} and Corollary \ref{C}, and justify Remark \ref{stable}(1). 

In preparation for the proofs of Theorem \ref{acm} and Corollary \ref{C} we note that, as $\mathbb F_2\times\mathbb F_2$ is not HS-stable by Theorem \ref{main}, there are sequences $U_{n,1},U_{n,2},V_{n,1},V_{n,2}\in\text{U}(d_n)$, for some $d_n\in\mathbb N$, such that 
\begin{enumerate}
\item $\|[U_{n,p},V_{n,q}]\|_2\rightarrow 0$, as $n\rightarrow\infty$, for every $1\leq p,q\leq 2$, and
\item $\inf_{n\in\mathbb N}(\|U_{n,1}-\widetilde U_{n,1}\|_2+\|U_{n,2}-\widetilde U_{n,2}\|_2+\|V_{n,1}-\widetilde V_{n,1}\|_2+\|V_{n,2}-\widetilde V_{n,2}\|_2)>0$, for any sequences $\widetilde U_{n,1},\widetilde U_{n,2},\widetilde V_{n,1},\widetilde V_{n,2}\in\text{U}(d_n)$ such that $[\widetilde U_{n,p},\widetilde V_{n,q}]=0$, for every $1\leq p,q\leq 2$.
\end{enumerate}
Consider the matricial ultraproduct $M=\prod_\omega\mathbb M_{d_n}(\mathbb C)$. Letting $U_p=(U_{n,p}),V_q=(V_{n,q})\in\mathcal U(M)$, condition (1) implies that $[U_p,V_q]=0$, for every $1\leq p,q\leq 2$.

\subsection{Proof of Theorem \ref{acm}} Assume by contradiction that the conclusion of Theorem \ref{acm} is false.

Let $1\leq p\leq 2$. Let $f:\mathbb T\rightarrow[-\frac{1}{2},\frac{1}{2}]$ is a Borel function satisfying $\exp(2\pi if(z))=z$, for all $z\in\mathbb T$, and define $h_p=f(U_p)$. 
 Then $h_p\in M$ is self-adjoint, generates the same von Neumann algebra as $U_p$ and satisfies $\|h_p\|\leq \frac{1}{2}$ and
 $U_p=\exp(2\pi ih_p)$. Similarly, let $k_p=f(V_p)$.

Let $A=h_1+ih_2$ and $B=k_1+ik_2$. Then $\|A\|\leq 1$ and $\|B\|\leq 1$. As $[h_p,k_q]=0$, for any $1\leq p,q\leq 2$, we have $[A,B]=[A,B^*]=0$. Represent $A=(A_n)$ and $B=(B_n)$, where $A_n,B_n\in\mathbb M_{d_n}(\mathbb C)$ satisfy $\|A_n\|,\|B_n\|\leq 1$, for every $n\in\mathbb N$. Then $\lim\limits_{n\rightarrow\omega}\|[A_n,B_n]\|_2=\lim\limits_{n\rightarrow\omega}\|[A_n,B_n^*]\|_2=0$. 

Since the conclusion of Theorem \ref{acm} is assumed false, we can find $A_n',B_n'\in\mathbb M_{d_n}(\mathbb C)$ such that $\lim\limits_{n\rightarrow\omega}\|A_n-A_n'\|_2=\lim\limits_{n\rightarrow\omega}\|B_n-B_n'\|_2=0$ and $[A_n',B_n']=[A_n,B_n'^*]=0$, for every $n\in\mathbb N$.
For $n\in\mathbb N$, denote by $P_n$ and $Q_n$ the von Neumann subalgebras of $\mathbb M_{d_n}(\mathbb C)$ generated by $A_n'$ and $B_n'$. Then $P_n$ and $Q_n$ commute and
 $\lim\limits_{n\rightarrow\omega}\|A_n-\text{E}_{P_n}(A_n)\|_2=\lim\limits_{n\rightarrow\omega}\|B_n-\text{E}_{Q_n}(B_n)\|_2=0$.

Then $\mathcal P=\prod_\omega P_n$ and $\mathcal Q=\prod_\omega Q_n$ commute, and $A\in \mathcal P,B\in \mathcal Q$. Hence $h_1,h_2\in \mathcal P$ and $k_1,k_2\in \mathcal Q$. 

Thus,  $U_p\in \mathcal U(\mathcal P)$ and $V_p\in \mathcal U(\mathcal Q)$, so we can find $\widetilde U_{n,p}\in\mathcal U(P_n)$ and $\widetilde V_{n,p}\in\mathcal U(Q_n)$, for every $n\in\mathbb N$, such that $U_p=(\widetilde U_{n,p})$ and $V_p=(\widetilde V_{n,p})$, for every $1\leq p\leq 2$. But then we have that $\lim\limits_{n\rightarrow\omega}\|U_{n,p}-\widetilde U_{n,p}\|_2=\lim\limits_{n\rightarrow\omega}\|V_{n,p}-\widetilde V_{n,p}\|_2=0$, for every $1\leq p\leq 2$. Since $[\widetilde U_{n,p},\widetilde V_{n,q}]=0$, for every $1\leq p,q\leq 2$, this contradicts (2), which finishes the proof. \hfill$\blacksquare$

\subsection{Proof of Corollary \ref{C}}
Let $P$ and $Q$ be the von Neumann subalgebras of $M$ generated by $\{U_1,U_2\}$ and $\{V_1,V_2\}$, respectively. Then $P$ and $Q$ commute. The last paragraph of the proof of Theorem \ref{acm} implies that there do not exist commuting von Neumann subalgebras $P_n,Q_n$ of $\mathbb M_{d_n}(\mathbb C)$, for all $n\in\mathbb N$, such that $P\subset\prod_\omega P_n$ and $Q\subset\prod_\omega Q_n$, which gives the conclusion.
\hfill$\blacksquare$

\subsection{Almost versus near commuting when one matrix is normal} The following result generalizes Remark \ref{stable}(1).

\begin{lemma}\label{ultra} Let $(M_n,\tau_n)$, $n\in\mathbb N$, be a sequence of tracial von Neumann algebras.
Let $x_n,y_n\in (M_n)_1$ such that  $y_n$ is normal, for every $n\in\mathbb N$, and  $\|[x_n,y_n]\|_2\rightarrow 0$. Then there are $x_n',y_n'\in M_n$ such that $x_n'y_n'=y_n'x_n'$ and $x_n'y_n'^*=y_n'^*x_n'$, for every $n\in\mathbb N$, and   $\|x_n-x_n'\|_2+\|y_n-y_n'\|_2\rightarrow 0$. 
 \end{lemma}
 
 This result can be proved qunatitatively by adapting \cite{Gl10,FK10}. Instead,  as in \cite{HL08}, we give a short proof using tracial ultraproducts.
 
 {\it Proof.} 
Consider the ultraproduct von Neumann algebra $M=\prod_{\omega}M_n$, where $\omega$ is a free ultrafilter on $\mathbb N$. Let $P$ and $Q$ be the von Neumann subalgebras of $M$ generated by $x=(x_n)$ and $y=(y_n)$.  Then $[x,y]=0$. Since $y$ is normal we get that $[x,y^*]=0$, so $P$ and $Q$ commute. Since $y$ is normal, we also get that $Q$ is abelian. By applying \cite[Theorem 2.7]{HS16} or \cite[Proposition C]{IS19} we can represent $x=(x_n')$ and $y=(y_n')$ so that the von Neumann subalgebras of $M_n$ generated by $x_n'$ and $y_n'$ commute, for all $n\in\mathbb N$. Since $\|x_n-x_n'\|_2+\|y_n-y_n'\|_2\rightarrow 0$, as $n\rightarrow\omega$,  the conclusion follows. 
 \hfill$\blacksquare$

\end{document}